\documentclass[12pt]{article}
\usepackage{amsmath, amsthm, amssymb}
\usepackage[pdftex]{graphicx}

\usepackage{mathrsfs}
\usepackage{color}
\def\Xint#1{\mathchoice
  {\XXint\displaystyle\textstyle{#1}}%
  {\XXint\textstyle\scriptstyle{#1}}%
  {\XXint\scriptstyle\scriptscriptstyle{#1}}%
  {\XXint\scriptscriptstyle\scriptscriptstyle{#1}}%
  \!\int}
\def\XXint#1#2#3{{\setbox0=\hbox{$#1{#2#3}{\int}$}
    \vcenter{\hbox{$#2#3$}}\kern-.5\wd0}}

\def\dashint{\Xint-}

\def\({\left(}
\def\){\right)}
\def\HH{{\mathcal{H}}}

\def\RR{{\mathbb{R}}}
\def\CC{{\mathbb{C}}}

\newcommand{\ZZ}{{\mathbb{Z}}}
\def\sm{\setminus}
\def\sm{\setminus}
\def\dm{\frac{1}{2}}

\def\MM{{\mathcal{M}_s}}

\newcommand{\weak}{\rightharpoonup}

\newcommand{\hev}{\vec{h}_{ex}}

\newcommand{\Om}{\Omega}

\newcommand{\ELD}{{\cal{L}}_{\eps,s}}
\newcommand{\ELDL}{{\cal{L}}_{\eps,s}^{\lambda}}

\newcommand{\eia}{e^{i\int_{z_{n-1}}^{z_n} A_z(x,y,z)\, dz}}

\newcommand{\eiaz}{e^{-i\int_z^{z_n} A_z(x,y,z')\, dz'}}

\newcommand{\grad}{\nabla}

\newcommand{\NN}{{\bf N}}

\newcommand{\eps}{\epsilon}

\renewcommand{\Im}{{\rm Im\,}}
\renewcommand{\Re}{{\rm Re\,}}

\newcommand{\he}{h_{ex}}

\newcommand\lep{|\ln \eps|}
\newcommand{\hb}{\overline{h}}
\newcommand{\Ab}{\overline{A}}

\newtheorem{thm}{Theorem}[section]
\newtheorem{lem}[thm]{Lemma}

\newtheorem{prop}[thm]{Proposition}
\newtheorem{rem}[thm]{Remark}
\newtheorem{defn}[thm]{Definition}

\newcommand{\be}{\begin{equation}}
\newcommand{\ee}{\end{equation}}
\newcommand{\bea}{\begin{eqnarray}}
\newcommand{\eea}{\end{eqnarray}}
\newcommand{\beann}{\begin{eqnarray*}}
\newcommand{\eeann}{\end{eqnarray*}}
\newcommand{\nnn}{\nonumber}

\begin{document}
\title{Minimizers of the Lawrence--Doniach Functional with Oblique Magnetic Fields}
\author{{\Large Stan Alama}\footnote{Dept. of Mathematics and Statistics,
McMaster Univ., Hamilton, Ontario, Canada L8S 4K1.  Supported
by an NSERC Research Grant.} \and  {\Large Lia Bronsard${}^*$}
  \and {\Large Etienne Sandier}\footnote{
Universit\'e Paris-Est,
LAMA -- CNRS UMR 8050, 
61, Avenue du G\'en\'eral de Gaulle, 94010 Cr\'eteil. France
\& Institut Universitaire de France.
{\tt sandier@u-pec.fr}
}}

\thispagestyle{empty}
\maketitle

\begin{abstract}

We study minimizers of the Lawrence--Doniach energy, which describes equilibrium states of superconductors with layered structure, assuming Floquet-periodic boundary conditions.  Specifically, we consider the effect of a constant magnetic field applied obliquely to the superconducting planes in the limit as both the layer spacing $s\to 0$ and the Ginzburg--Landau parameter $\kappa=\eps^{-1}\to\infty$, under the hypotheses that $s=\eps^\alpha$ with $0<\alpha<1$.  By deriving sharp matching upper and lower bounds on the energy of minimizers, we determine the lower critical field and the orientation of the flux lattice, to leading order in the parameter $\eps$.  To leading order, the induced field is characterized by a convex minimization problem in $\RR^3$.  We observe a ``flux lock-in transition'', in which flux lines are pinned to the horizontal direction for applied fields of small inclination, and which is not present in minimizers of the anisotropic Ginzburg--Landau model.  The energy profile we obtain suggests the presence of ``staircase vortices'', which have been described qualitatively in the physics literature.

\bigskip

\noindent
{\bf Keywords:} Calculus of variations, elliptic equations and systems, superconductivity, vortices.

\medskip

\noindent
{\bf MSC subject classification:}  35J50, 58J37

\end{abstract}

\newpage


\baselineskip=18pt

\section{Introduction}

One characteristic of the high-temperature superconductors is their strong anisotropy, which distinguishes them from conventional superconducting materials.  Indeed, the cuprate superconductors (such as BSCCO, Ba${}_2$Sr${}_2$Ca${}_n$Cu${}_{n+1}$O${}_{2n+6+x}$, or TBCCO, Tl${}_2$Ba${}_2$Ca${}_2$Cu${}_3$O${}_{10}$) are crystalline materials with a distinct layered structure.  Conduction of superconducting electrons is favored in cuprate planes, which are separated by insulating layers in the crystal structure.  Physicists recognize that the most appropriate description of such layered superconducting structures is given by the the Lawrence--Doniach model \cite{LaDo}, which was introduced in 1971 (well before the advent of the high-$T_C$ materials.)  The Lawrence--Doniach energy treats the superconducting solid as an array of two-dimensional superconducting planes described by a Ginzburg--Landau energy, with discrete coupling of adjacent layers in the orthogonal direction.  In the limit as the spacing between adjacent layers tends to zero, the Lawrence--Doniach energy converges to the anisotropic three-dimensional Ginzburg--Landau energy \cite{CDG,BK}.  This approximation is adequate for certain high-$T_C$ materials, but not for the more anisotropic varieties \cite{Iye}, where it is believed that the layered structure plays a strong role.  

Many important questions in superconductivity concern the response of a superconducting sample to an applied magnetic field.  In particular, superconductors expell an applied magnetic field of low strength (the Meissner effect,) and there is a critical value of the applied field's magnitude $H_{C1}$ at which flux lines penetrate the material via vortices.  In analogy with our previous paper \cite{ABS2} on the three-dimensional anisotropic Ginzburg--Landau model, we pose the following question:  for a given orientation of the applied magnetic field (with respect to the superconducting planes,) what is the value of $H_{C1}$, and in which direction do vortex lines and lines of magnetic flux lie?  As in \cite{ABS2}, we choose a triply-periodic setting in $\RR^3$.  This represents a superconductor occupying all of space, and periodic vortex lattices (as are normally observed in the interior of actual superconducting samples.)  

We now define the analytical setting for our problem.  The superconducting planes $\mathcal{P}_n$, $n\in\ZZ$, are assumed to be horizontal, with uniform spacing $s>0$,
$$   \mathcal{P}_n = \RR^2 \times \{z=z_n:=ns\}, \qquad n\in \ZZ. $$
On each plane is defined a complex-valued order parameter $u_n: \ \mathcal{P}_n\to \CC$.  As in the Ginzburg--Landau model, $|u_n|=1$ represents the purely superconducting state, and $u_n=0$ the normal state. 
The superconducting currents interact with a magnetic field, described by a vector potential $A: \RR^3\to \RR^3$ via $h=\nabla \times A$.  

We assume periodicity with respect to a fundamental domain $\Omega\subset\RR^3$ defined by a fixed basis $\{\vec{v_1}, \vec{v_2}, \vec{v_3}=(0,0,L)\}$ with $\vec{v_3}\perp \vec{v_1},\vec{v_2}$.  The horizontal period is arbitrary, while the vertical period is dictated by the discrete geometry of the system.  The geometry of the period cell dictates that the distance $s$ between adjacent planes be compatible with the period, so we assume
$$   s = {L\over N} \quad\text{for $N\in \mathbb{N}$}.  $$
We also define the restriction of each plane to $\Omega$,
$$  P_n = \mathcal{P}_n \cap\Omega 
  = \{  t_1\vec{v_1}+ t_2\vec{v_2} + (0,0,z_n): \ t_1,t_2\in [0,1] \}, \quad  n=1,\dots,N.  $$
As our analysis will show, in our asymptotic regime the natural period of energy minimizers will be much smaller, and the choice of fixed period cell $\Omega$ will not enter into the geometry of the minimizers at all.

To write the LD energy we introduce some convenient notation.  We denote by $\grad'= (\partial_x,\partial_y)$, and $A'=(A_x(x,y,z), A_y(x,y,z))$, 
$A'_n=(A_x(x,y,z_n), A_y(x,y,z_n))$.  We write $(u_n, A)$ as a shorthand for $(\{u_n\}_{n\in\ZZ}, A)$.
Then, the energy in the period domain $\Omega$ may be written as:
\bea\label{LD}
\ELDL (u_n, A)&=&
    s\sum_{n=1}^{N(s)}
     \int_{P_n} \left[  \frac12
            \left|\left( \grad' -i 
A'_n\right)u_n  \right|^2
+ \frac{1}{4\eps^2} (|u_n|^2-1)^2  \right]\, dx\, dy
          \\
  \nnn      &&\qquad + s\sum_{n=1}^N
                    \int_{P_n}  {1\over 2\lambda^2 s^2}\left| u_n -
                      u_{n-1}\eia
                           \right|^2
                   \, dx\, dy \\
\nnn       &&
             \qquad   +\
                \frac12 \int_\Om \left| \nabla\times A- \hev\right|^2 \, dx\, dy\, dz .
   \eea
Here, $\eps>0$ represents the reciprocal of the Ginzburg--Landau parameter; we will assume $\eps\ll 1$, which is typical for type-II superconducting materials.  The constant $\lambda>0$ represents the Josephson penetration depth, and will assumed to be fixed in this paper.  (Results on the large $\lambda$ limit of the LD minimizers may be found in \cite{ABB1, ABB2}.)  The applied field $\hev$ is a given constant vector, which will  depend on $s$ and $\eps$.

Next we must define a space of functions for $\ELDL$.  
\begin{defn}\label{deffloq}
  We say $(u_n,A)\in \HH$ if $u_n\in H^1_{loc}(\RR^2;\CC)$ for all $n\in\ZZ$, $A\in H^1_{loc}(\RR^3;\RR^3)$, and
there exist functions $\omega_j\in H^2_{loc}(\RR^3)$, $j=1,2,3$ so that:
\be\label{tH}
\left\{  \begin{gathered}
  u_n(\vec{x}+\vec{v}_j)= u_n(\vec{x})e^{i\omega_j(\vec x,z_n)}, \quad j=1,2, \ n\in\ZZ;  \\
  u_{n+N}(\vec{x}) = u_n(\vec{x})e^{i\omega_3(\vec{x},z_n)}, \quad n\in\ZZ; \\
  A(\vec{x}+\vec{v}_j)= A(\vec{x}) + \nabla\omega_j(\vec{x}), \quad j=1,2,3
\end{gathered} \right.
\ee
holds for all $\vec{x}=(x,y,z)=(x',z)\in \RR^3$. 
\end{defn}
 That is, the configuration $(u_n,A)\in\HH$ is $\Om$-periodic up to gauge transformation.  In particular, the gauge-invariant quantities such as the energy density and the currents, density, and induced magnetic field,
$$
\begin{gathered}   j'_n(x,y)=\Im\{\overline{u_n}\,(\grad'-iA_n')u_n\},  \quad
    j_z(x,y)=-\Im\left\{ \overline{u_n} \, u_{n-1}\eia\right\} \\
    \rho_n(x,y)=|u_n|, \quad h=\nabla\times A,  
\end{gathered}    $$
are all $\Om$-periodic.

An important consequence of periodic boundary conditions is the exact quantization of the magnetic flux through period planes:  see Proposition~\ref{fluxprop} below.  Indeed, our analysis of minimizers for the Lawrence--Doniach model is dependent on a sharp lower bound derived using the flux quantization and a slicing method.

\medskip

The high-$T_C$ materials being strongly type-II superconductors, it is natural to consider the London limit $\eps\to 0$.  In addition, the distance between layers being very small, we will also consider $s\to 0$.  In addition, we will assume throughout the paper that
\begin{equation}\label{seps}   s = \eps^\alpha  
\end{equation}
for some $\alpha\in (0,1)$.  This choice is interesting for two reasons.  First, $\eps\ll s$ is needed for the discrete nature of the Lawrence--Doniach model to be retained in the limit.  Indeed, if $s=O(\eps)$ our result predicts the same leading order behavior for minimizers of Lawrence--Doniach as for the anisotropic Ginzburg--Landau model, obtained in \cite{ABS2}.  Secondly, our results in \cite{ABS3} show that for applied fields imposed orthogonally to the superconducting planes, the lower critical field (the smallest field strength at which flux lines penetrate the sample) is of order $\lep$, whereas for applied fields parallel to the layers, the lower critical field is of order $|\ln s|$.  With the choice \eqref{seps}, both critical fields are of the same order of magnitude, and so we expect both components of an oblique applied field to contribute to the energy at the same scale.  Thus, we assume the external applied fields are of the form
\begin{equation}\label{hex}   h^\eps_{ex} = H_{ex}\,\lep
    =\alpha^{-1} H_{ex}\, |\ln s|,  
\end{equation}
with $H_{ex}$ a fixed, constant vector (or 2-form) in $\RR^3$.  In fact, it will be convenient to think of vectors $H\in\RR^3$ as 2-forms, with the association 
$$(H_1,H_2,H_3) \leadsto H_1 dx^2 dx^3 + H_2 dx^3 dx^1 + H_3 dx^1 dx^2.  $$
Employing a Riemannian metric $g=\text{diag}\, (1, 1, \lambda^2)$ to represent the anisotropy in the directions in the Lawrence--Doniach model, we define norms on the 2-form $H$,
\begin{gather*}   \|H\|_g = \sup_{X,Y\in\RR^3} {H(X,Y)\over |X|_g |Y|_g} =
                     \lambda^{-1} \sqrt{ H_1^2 + H_2^2 + \lambda^2 H_3^2},
   \\
     \|H\|_{g^{-1}}
        = \sup_{X,Y\in\RR^3} {H(X,Y)\over |X|_{g^{-1}} |Y|_{g^{-1}}} =
                     \lambda \sqrt{ H_1^2 + H_2^2 + \lambda^{-2} H_3^2}.
\end{gather*}

Our main result is the following:
\begin{thm}\label{bigthm}
Assume \eqref{hex}, and \eqref{seps}.  Let $(u_{n,\eps},A_\eps)\in \HH$ be minimizers of $\ELDL$.  Then, 
$$   {h_\eps\over\lep}\weak H_* \quad \text{in $L^2(\Omega)$},  $$
with $H_*\in \RR^3$ constant, and
\begin{equation}\label{conv}   \lim_{\eps\to 0} {\ELDL(u_{n,\eps},A_\eps)\over |\Omega| \, \lep^2} = \frac12 \left( (1-\alpha) |H_*\cdot e_3| +
    \alpha \| H_*\|_g  + |H_*-H_{ex}|^2\right) =: F(H_*).
\end{equation}
Moreover, $H_*$ minimizes $F(H)$ among vectors in $\RR^3$.
\end{thm}

In Proposition~\ref{dual}, we use Fenchel duality to obtain a more geometrical characterization of the minimizing field $H_*$.  We show that the optimal $H_*$ is in fact the closest point to the origin in the closed convex body $K+H_{ex}$ where $K$ is defined by
$$   U\in K \iff \begin{cases}
\| U-{1-\alpha\over 2} e_3\|_{g^{-1}}\le {\alpha\over 2}, 
  &\text{if $U_3\ge {1-\alpha\over 2},$} \\
\| U+{1-\alpha\over 2} e_3\|_{g^{-1}}\le {\alpha\over 2}, 
  &\text{if $U_3\le -\left({1-\alpha\over 2}\right),$} \\
\| U'\|_{g^{-1}}\le {\alpha\over 2} ,
   &\text{if $-\left({1-\alpha\over 2}\right) \le U_3\le {1-\alpha\over 2},$}
\end{cases}
$$
\begin{figure}[htbp] 
   \centering
   \includegraphics[width=3.5in]{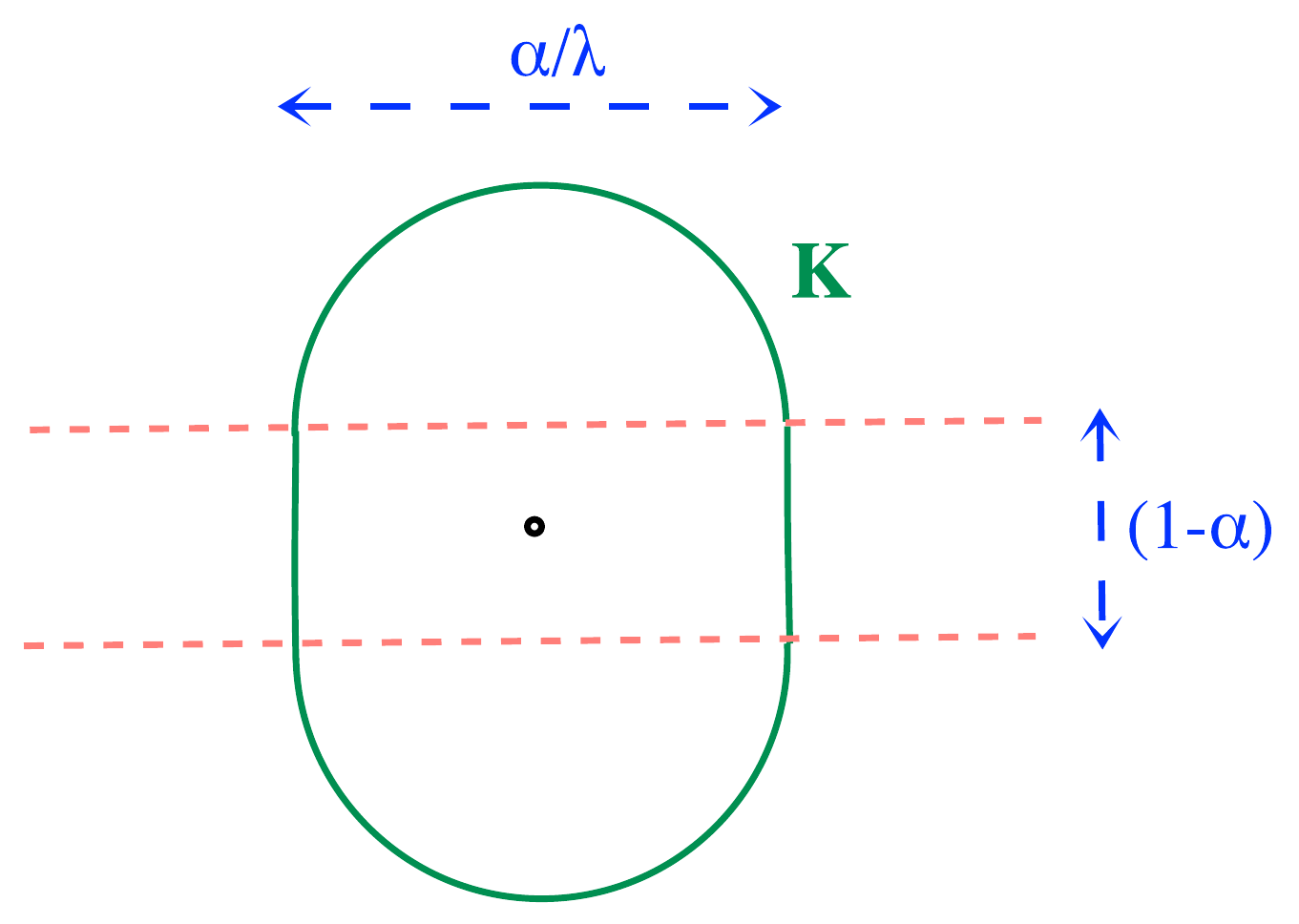} 
   \caption{The convex set $K$.  }
   \label{fig1}
\end{figure}

That is, the set $K$ has a cylindrical side of height $(1-\alpha)$, and ellipsoidal caps (see figure~\ref{fig1}.)  The presence of the cylindrical sides of $K$ provides for an interesting phenomenon in the Lawrence--Doniach model which is not observed for minimizers of the anisotropic Ginzburg--Landau energy.  Indeed, let us consider applied fields 
$H_{ex}$ with a fixed angle $\theta$ with respect to the horizontal superconducting planes, 
$$H_{ex}= (H'_{ex}, H_{ex,3})= |H_{ex}|\cos\theta\,\vec{e_1} + |H_{ex}|\sin\theta\, \vec{e_3},  $$
allowing the magnitude $|H_{ex}|$ to vary.

For sufficiently small $\theta$, $|\tan\theta|<\lambda\left({1-\alpha\over\alpha}\right)$ to be precise, the point $H$ of $K+H_{ex}$ nearest to the origin can occur in three different locations, each with its own interpretation.  The situation is illustrated in figure~\ref{fig2} below.
For $|H'_{ex}|=|H_{ex}\cos\theta|\le {\alpha\over 2\lambda}$, $H=0$ lies in the set $K+H_{ex}$ itself (left part of figure~\ref{fig2}.)  This represents the Meissner effect, and the magnetic field is expelled from the sample (no vortices.)

  For ${\alpha\over 2\lambda|\cos\theta|}<|H_{ex}|\le {1-\alpha\over 2\lambda|\sin\theta|},$ the nearest point in $K+H_{ex}$ to the origin lies along the vertical edge of the set, and so $H$ is  horizontal (center image of figure~\ref{fig2}.)  Thus, although the applied field is not horizontal, flux lines are trapped by the superconducting planes and the vertical component is suppressed.  This phenomenon has been observed in the layered high-$T_C$ superconductors \cite{SKMMW}, and is known as a ``flux lock-in transition.''  It is an example of a physical phenomenon which is not present in the anisotropic Ginzburg--Landau model, and used as a justification for the Lawrence--Doniach model in the physics literature.  We note that recent work by Bauman \& Xie \cite{BaX} confirms the flux lock-in phenomenon in a different asymptotic regime, as $\lambda\to\infty$ with $s,\eps$ fixed (as in \cite{ABB2}), in nearly-parallel applied fields.

  For $|H_{ex}| > {1-\alpha\over 2\lambda|\sin\theta|}$, the nearest point lies on the ellipsoidal lower part of the set $K+H_{ex}$, and so $H$ gradually lifts away from the horizontal direction.  (See the right image of figure~\ref{fig2}.)  We note that for larger angles, $|\tan\theta|\ge\lambda\left({1-\alpha\over\alpha}\right)$, the closest point of $K+H_{ex}$ will always lie either inside the set or on one of the ellipsoidal caps, and so there is no flux pinning phenomenon observed in that case.

\begin{figure}[htbp] 
   \centering
   \includegraphics[width=5.5in]{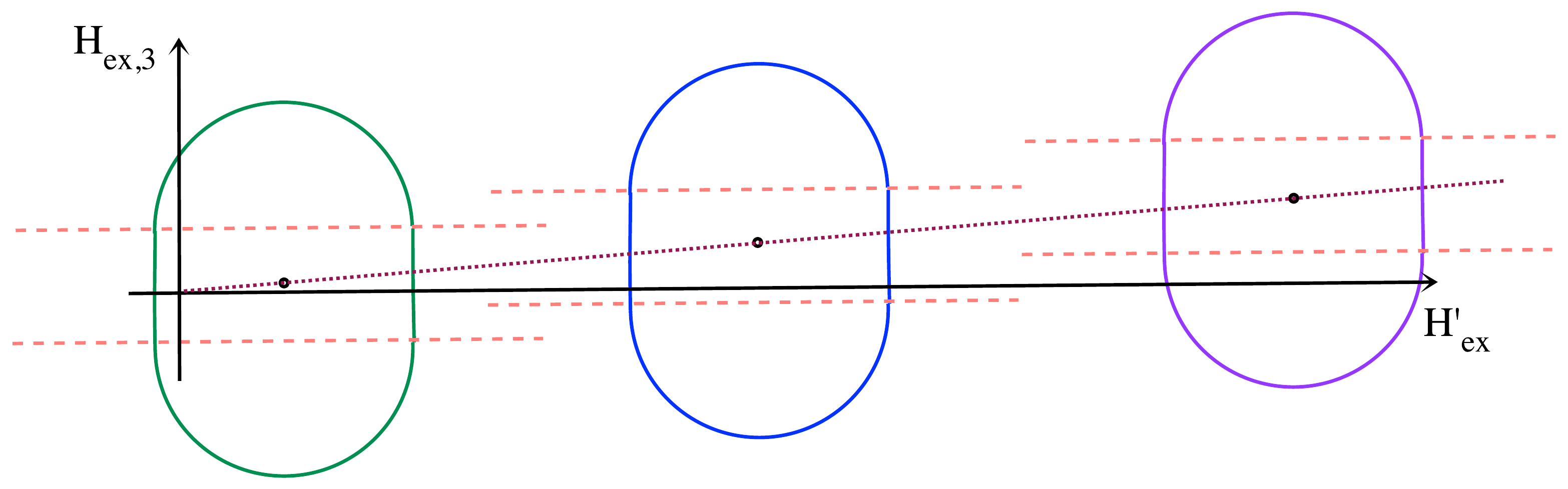} 
   \caption{Determining the induced field for external fields $H_{ex}$ applied at a small fixed angle $\theta$, with $|\tan\theta|<\lambda\left({1-\alpha\over\alpha}\right)$ with respect to the horizontal.  }
   \label{fig2}
\end{figure}  
  
From the form of the set $K$ we may thus obtain a formula for the lower critical field $H_{c1}(\theta)$, the largest value of $|H_{ex}|$ for which $H=0$ is the induced field of minimizers:
$$   H_{c1}(\theta) = \begin{cases}
   {\alpha\over 2\lambda\cos\theta}, 
       & \text{if 
        $\tan\theta\le\lambda\left({1-\alpha\over\alpha}\right)$,} 
         \\
    {(1-\alpha)\sin\theta 
       + \sqrt{\alpha^2\sin^2\theta - (1-2\alpha)\lambda^2\cos^2\theta}
         \over 2(\lambda^2\cos^2\theta + \sin^2\theta)},
         &\text{if
           $\tan\theta>\lambda\left({1-\alpha\over\alpha}\right)$,}
\end{cases}
$$
for $\theta\in [0,{\pi\over 2}]$.
  
\medskip
  
It is natural to ask what shape the vortex lines will have in this regime.  
In the periodic case, and with applied fields on the order $\lep$, the result we obtain captures the average field over a dense vortex lattice, and so a detailed characterization of the magnetic flux lines is not available to us.  Moreover, since the order parameters $u_n$ are only defined on the superconducting planes, the usual description of vortices as zeros of a continuous order parameter is not applicable here.  Nevertheless, an image of the rough structure of minimizers may be deduced from our upper and lower bound constructions and the form of the leading order term in the energy.  
First, since $s\to 0$, to the first approximation the energy resembles the anisotropic Ginzburg--Landau energy, except that the relevant length scale is $s$, determined by the coupling term in the energy.  Thus, each flux line contributes an energy of the order $\pi|\ln s|=\alpha\pi\lep$ per unit length.  The second term in the limiting energy $F(H)$ (see \eqref{conv})) arises in this way.  However, when a flux line penetrates one of the superconducting planes, it leads to a greater energy cost.  Inside each vortex tube of width $O(s)$, we observe a two-dimensional Ginzburg--Landau vortex with scale $\eps$, yielding an energy of $\pi\ln(s/\eps)=\pi(1-\alpha)\lep$ per vortex.  This is the origin of the first term in $F(H)$.  This leads us to imagine that the flux lines form a step-like pattern: lying nearly horizontally between planes, where the energy per length of each flux line is smaller; then penetrating each plane nearly vertically, to produce two-dimensional vortices in each superconducting layer.  Thus we obtain a first mathematical confirmation of the appearance of ``staircase vortices'', which were described qualitatively in the physics literature \cite{Bulk, Iye}, as illustrated in figure~\ref{fig3} below.

\begin{figure}[htbp] 
   \centering
   \includegraphics[width=2.5in]{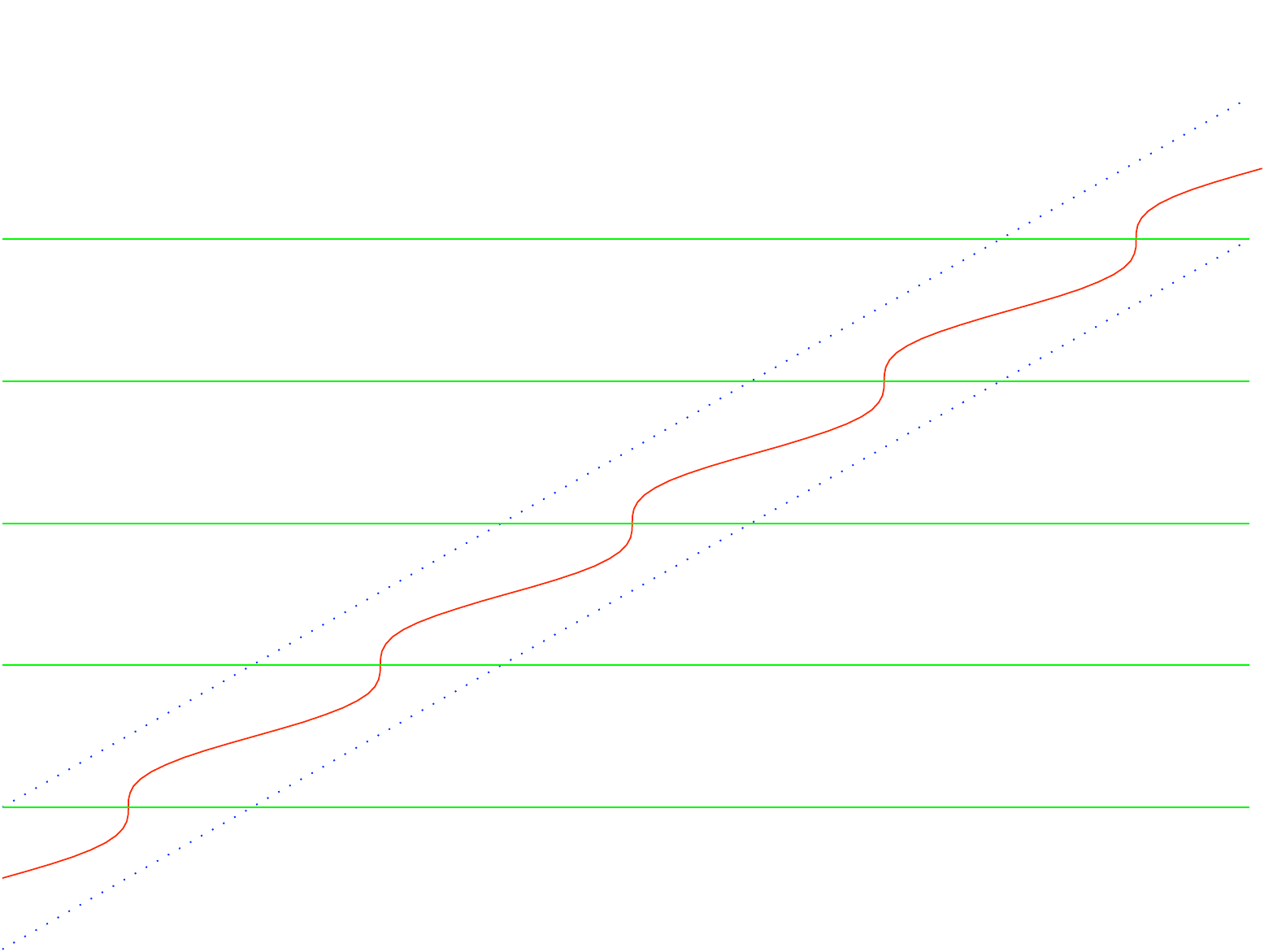} 
   \caption{A staircase vortex:  the flux line passes through each plane orthogonally, becoming flat in each gap to approximate the limiting direction $H_*$.}
   \label{fig3}
\end{figure}

The Lawrence--Doniach model contains several dimensionless parameters, and so there are many interesting asymptotic regimes to be considered.  Our previous papers \cite{ABS1, ABS3} consider applied fields which are parallel to the principal axes, in a simultaneous limit $s, \eps\to 0$.  Applied fields parallel to the superconducting planes are particularly interesting, and we identify two very different limiting regimes.  The first regime,  $s^2 |H_{ex}|\ll 1$, is Ginzburg--Landau-like, with a regular array of vortices whose mutual distance $\delta\sim |H_{ex}|^{-1/2}\gg s$, the radius of a vortex core.  The second asymptotic regime is $s^2 |H_{ex}|\gg 1$, for which the sample appears transparent to the penetration of the applied field with no loss of superconductivity, and many flux lines pass through each gap.  Our hypothesis $|H_{ex}|=O(|\ln \eps|)=O(|\ln s|)$ clearly belongs to the first case, and the result of Theorem~\ref{bigthm} is consistent with the Ginzburg--Landau-like image of minimizers, although with a more subtle geometry.  The leading order behavior of minimizers in other asymptotic regimes, both for oblique and parallel applied fields, remains open.

The paper is organized as follows:  section 2 contains some preliminary definitions and results concerning the variational formulation of the periodic problem, and the use of differential forms.  In section 3 we derive an upper bound on the energy of minimizers by modifying a basic construction for the anisotropic Ginzburg--Landau energy in \cite{ABS2}.  A complementary lower bound is obtained in section 4, using flux quantization and a slicing method.  Section 5 is devoted to showing the equivalent geometrical formulation of the limiting problem for the normalized induced field, $H$.

\section{Preliminary results}
\setcounter{thm}{0}


As mentioned above, given the anisotropy in the model it will be convenient in most of the paper to express various quantities as differential forms rather than vector fields, for ease in applying transformations to evaluate integrals.  To our vector potential, we associate a differential 1-form which we also denote by $A$,
$$  A= A(x,y,z) = A_x \, dx + A_y\, dy + A_z\, dz.  $$
The induced magnetic field is then a 2-form, $h= dA$.  As above, to express the Lawrence--Doniach energy we also define the notation
$$   d'= \partial_x\, dx + \partial_y\, dy, \qquad
    A'_n = A'_n(x,y) = A_x(x,y,z_n)\, dx + A_y(x,y,z_n)\, dy.  $$
The energy functional may then be rewritten in an equivalent way as:
\begin{align*}
\ELDL (u_n, A)&=
    s\sum_{n=1}^{N(s)}
     \int_{P_n} \left[  \frac12
            \left| d'u_n -i 
A'_nu_n  \right|^2
+ \frac{1}{4\eps^2} (|u_n|^2-1)^2  \right]\, dx\, dy
          \\
     &\qquad + s\sum_{n=1}^N
                    \int_{P_n}  {1\over 2\lambda^2 s^2}\left| u_n -
                      u_{n-1}\eia
                           \right|^2
                   \, dx\, dy \\
       &
             \qquad   +\
                \frac12 \int_\Om \left| d A- h_{ex}^\eps\right|^2 \, dx\, dy\, dz .
\end{align*}  
We continue to denote by $(u_n,A)\in\HH$ (as in Definition~\ref{deffloq}) the class of admissible Floquet-periodic configurations for $\ELDL$, expressed with $A$ as a 1-form and the differential $d\omega_j$ replacing the gradient $\nabla\omega_j$ (j=1,2,3) in \eqref{tH}.

In order to assert the existence of energy minimizers we must fix a gauge in which the Lawrence--Doniach energy is coercive.  
For a {\em constant} 2-form,
$ \hb = \frac12\sum_{i,j,k} \eps_{ijk}\,\hb_i\, dx^j\wedge dx^k$, we select as representative,
\be\label{Abar}
\Ab = \frac12 \sum_{i,j,k}\, \eps_{ijk} \,\hb_i\, x_j\, dx^k.
\ee
That is, if we think of
$\hb = (\hb_1, \hb_2, \hb_3)$ as a vector field, then $\Ab = \frac12\, \hb\times \vec{x}$.  
We then say $(u_n,A)\in \HH_*$ if $u_n\in H^1_{loc}(\RR^2;\CC)$ for all $n\in\ZZ$, $A\in H^1_{loc}(\RR^3;\Lambda^1(\RR^3))$, and
there exist a {\em constant} 2-form $\hb\in\Lambda^2(\RR^3)$ such that:
\be\label{tH*}  \left.
\begin{gathered}
A= \Ab + A_0, \quad d^* A_0=0,
\quad A_0(\vec{x} +\vec{v}_j)=A_0(\vec{x}), \quad j=1,2,3, \ \vec{x}\in\RR^3,   \\
u_n(\vec{x}+\vec{v}_j) = u_n(\vec{x})e^{-i\Ab(\vec{x},z_n)\cdot\vec{v}_j}, \qquad j=1,2, \ \vec{x}\in\RR^2, \ n\in\ZZ \\
u_{n+N}(\vec{x}) = u_n(\vec{x})e^{-i\Ab(\vec{x},z_n)\cdot\vec{v}_3}, \quad \vec{x}\in\RR^2, \ n\in\ZZ,
\end{gathered}  \right\}
\ee
where $\Ab$ is associated to $\hb$ as in (\ref{Abar}).
Note that $(u_n,A)\in \HH_*$ satisfy (\ref{tH}) with $\omega_j(\vec{x})= -\Ab(\vec{x})(\vec{v}_j)$.  
\begin{lem}\label{gauge1}
\begin{enumerate}
\item[{\bf (a)}]  For any $(u_n,A)\in\HH$ there exists $\gamma\in H^2_{loc}(\RR^3)$ so that
$$   (u_n e^{i\gamma(\cdot, z_n)}, \ A+d\gamma)\in \HH_*.  $$
Moreover, the constant form $\hb$ is the componentwise average value of $d A$ in $\Om$.
\item[{\bf (b)}]  There exists a constant $C_0=C_0(\Om)$ such that for any $(u_n,A)\in \HH_*$
\be\label{Abound}   \| A\|_{H^1(\Om)} \le C_0 \| d A\|_{L^2(\Om)}.   
\ee
\end{enumerate}
\end{lem}

This has been proven in our previous paper \cite{ABS3}.

An important consequence of Floquet periodic boundary conditions is the quantization of flux:

\begin{prop}\label{fluxprop}
Let $(u_n, A)\in \HH$ and assume that for some $n$, $|u_n|>0$ on $\partial P_n$.  Then
$$   \deg\left( {u_n\over |u_n|};\partial P_n\right) 
   ={1\over 2\pi} \int_{P_n} dA.
$$
\end{prop}

Given the estimate (\ref{Abound}), the Lawrence--Doniach functional is coercive on $\HH_*$ and the existence of minimizers for $\ELDL$ in $\HH_*$, for fixed values of the parameters $s=L/N$, $\lambda$, $\eps$, $h^\eps_{ex}$, follows from standard arguments \cite{Aydi,Odeh}.  Minimizers satisfy a system of Euler--Lagrange equations which contain singular terms along the planes $P_n$.  The regularity of the order parameters $u_n$ and partial regularity of the magnetic potential $A$ (in the Coulomb gauge) has been studied by Bauman \& Ko \cite{BK}.  Our treatment will be largely variational, and we will only require the following estimate on minimizers, which is proven via the maximum principle (see \cite{BK}):

\begin{prop}\label{apriori}
Let $(u_n,A)\in\HH_*$ be critical points of $\ELDL$.  Then $|u_n|\le 1$ for all $x\in P_n$, for all $n\in\NN$.
\end{prop}

\section{Upper bound}\label{UBsection}
\setcounter{thm}{0}

In this section we assume that 
\begin{equation}\label{truc} s=\eps^\alpha,\quad h_{ex}^\eps = H_{ex}\lep = {1\over\alpha} H_{ex}|\ln s|,  \end{equation}
where $H_{ex}$ is a constant  2-form, and $0<\alpha<1$. We have 

\begin{prop}\label{ubprop} Assuming \eqref{truc}, for any constant 2-form $h\in \Lambda^2(\RR^3)$ and for any $s>0$, there exists 
 $(u_n^s, A^s)\in \HH$ such that
$$   {\ELDL(u_n^s, A^s)\over |\Omega| |\ln\eps|^2}
           \le \frac12\left[ (1-\alpha) |h_3| 
           + \alpha\| h\|_g + | h - H_{ex}|^2\right]
$$
\end{prop}

\begin{rem}  The case $h_3:= h(\vec{e_1},\vec{e_2})=0$ --- which corresponds to the case of  a magnetic field parallel to the superconducting planes --- was treated in Proposition~4.4 of \cite{ABS3}, thus  we assume throughout the proof that 
\be\label{h3}
    h_3= h(\vec{e_1},\vec{e_2})\neq 0  
\ee
\end{rem}

\begin{proof}[Step~1, the reference configuration] 
 We recall the construction from Proposition~3.1 of \cite{ABS2}, but with $s$, $H_{ex}/\alpha$ and $h/\alpha$  playing the respective roles of $\eps$, $h_{ex}$, $h$. We begin by defining a reference configuration $(v,B)$ in $\RR^3\setminus (\ZZ^2\times\RR)$.
 
Let $K=\left[ -\frac12,\frac12\right]^2$ and $f$ solve
$$\left\{  \begin{array}{ll}
                       -\Delta f = 2\pi(\delta_{(0,0)} -1), & \mbox{in $K$,} \\
                     \  \partial_\nu f =0 & \mbox{on $\partial K$,}
                \end{array}  \right.
$$
We let 
$$   B=2\pi\, x\, dy,  $$
so that  $dB = 2\pi dx\wedge dy$. 
Then we define $v$ by
\begin{equation}\label{defphi}  \text{$v= e^{i\phi}$, where $d\phi = j+ B$ and $j=*df$,}\end{equation} 
where $*dx = dy$ and $*dy=-dx$. Given $p\in\ZZ^2$ we have $f(x) = \log|x-p| + f_p(x)$, where $f_p$ is smooth in a neighbourhood of $p$ and $\Delta f_p = 2\pi$, hence $j = d\theta_p + *df_p$, where $\theta_p$ is the angle in polar coordinates centered at $p$. It follows that $j+B = d\theta_p + \omega_p$, where $d\omega_p = 0$ and $\omega_p$ is smooth near $p$. Thus $v=e^{i\phi}$ is well defined and for each $p\in\ZZ^2$  we have 
\begin{equation}\label{vest1}
v=e^{i(\theta_p + \xi_p)} 
\end{equation}
where $\xi_p$ is smooth near $p$. 

Following the proof of Proposition~3.1 of \cite{ABS2}, we 
define $(v_s,B_s)$.  Let $\rho_s$ be a smooth cutoff first defined in the square $K$  equal to $1$ outside $B(0,2s)$ and equal to $0$ in $B(0,s)$, and then extended by periodicity to $\RR^2$. Then extend  the definition of $\rho_s$ and $v$ to $\RR^3$ by making them $z$-invariant and let 
\begin{equation}\label{aniso}
  v_s=(\rho_s v)\circ\widehat\Phi_s, \qquad B_s= \widehat\Phi^*_s B.  
\end{equation}
Here  $\widehat\Phi_s = |\ln s|^{1/2}\Phi_0$, where $\Phi_0$ is linear, independent of $s$ and defined by 
$$\Phi_0(\vec b_i) = \vec e_i, $$
with $(\vec b_i)_{i=1,2,3}$ a $g$-orthogonal basis of $\RR^3$ such that $\vec b_3$ is in the kernel of $h$ --- i.e. such that $h(\vec{b_3},\vec{c})=0$ for all $\vec{c}$ --- and normalized by the requirement $ \frac h\alpha(\vec b_1, \vec b_2)=2\pi$.  More precisely 
\bea\nnn
|\vec b_i|_g = |\vec b_j|_g \quad\forall i,j, \qquad g(\vec b_i,\vec b_j)=0 \quad i\neq j ,\\
\label{hb3}
h(\vec b_3,\cdot)\equiv 0, \qquad h(\vec b_1, \vec b_2)=2\pi\alpha. 
\eea
Note that by \eqref{h3} the vector $\vec b_3$ cannot lie in the horizontal plane, hence 
\bea\label{vertical}\Phi_0^{-1}(\vec e_3)\notin \text{Span}(\vec v_1,\vec v_2).\eea

We have  $dB_s = \widehat\Phi^*_s (dB) = \widehat\Phi^*_s (2\pi\,dx\wedge dy) = |\ln s| h/\alpha$, and the anisotropic Ginzburg-Landau energy of the configuration in a period cell
$$  \widehat C_s = \widehat\Phi_s^{-1}\left( \left[-\frac12,\frac12\right]^3\right),   $$
following the computations in equations (14)--(16) of \cite{ABS2},  satisfies
\be\label{vs}  {G_s(v_s, B_s; \widehat C_s) \over |\widehat C_s| \, |\ln s|^2} \le \frac12 \left(
\left\| {h\over\alpha}\right\|_g + \left| {h\over\alpha}- {H_{ex}\over\alpha}\right|^2\right) + o(1), 
\ee
where 
$$   G_s(v, B;\Omega) = \int_{\Omega}
\left\{
 \frac12 |(d - i B)v|_g^2 + {1\over 4 s^2}(|v|^2-1)^2 +
      \frac12 \left| dB - h^\eps_{ex} \right|^2\right\}  dV_e.
$$
Note that we  integrate all terms with respect to the euclidean volume element  $dV_e$. 
\end{proof}

\begin{proof}[Step~2, Periodicity] 
The configuration $(v_s,B_s)$ can be perturbed to be $\mathcal{A}$-periodic: We define the set of $\mathcal A$-periodic configurations $\HH$ by saying  $(v,B)\in \HH$ if 
$(v, B)\in H^1_{loc}(\RR^3;\mathbb{C})\times 
   H^1_{loc}(\RR^3;\Lambda^1(\RR^3))$
and there exist
$\omega_j\in H^2_{loc}(\RR^3)$, $j=1,2,3$ for which 
$$  v(\cdot + \vec{v_j}) = v(\cdot)\exp (i\omega_j), \qquad
    B(\cdot + \vec{v_j}) = B(\cdot)+ d\omega_j, \quad j=1,2,3.
$$

Then $(v_s,B_s)$ is made $\mathcal A$-periodic by using rational approximation of the transformation matrices.  We know that  $(v_s,B_s)$ is periodic w.r.t. $\mathcal B_s=(\vec b_i/|\ln s|^{1/2})$. Then we let $\mathcal B'_s$ be the basis such that
$$M_{\mathcal B'_s\mathcal  A}= [M_{\mathcal B_s \mathcal A}] = [|\ln s|^\dm M_{\mathcal B \mathcal A}],  $$
where brackets denote the integer part of each matrix entry, where $M_{\mathcal A\mathcal B}$ is the matrix whose columns are the coordinates of the vectors of $\mathcal B$ in the basis $\mathcal A$, and where we have set $\mathcal B = (\vec b_i)_i$. We define the configuration $(u_s, A_s)$ by 
$$   u_s = (\rho_s v)\circ \Phi_s, 
   \quad  A_s = \Phi_s^* B,  $$
where $\Phi_s$ is the linear map which maps the $i$-th vector in $\mathcal B'_s$ to $|\ln s|^{1/2} \vec e_i$. It is a perturbation of $\widehat\Phi_s$ in that $\Phi_s = \widehat\Phi_s\circ\Psi_s$, where  $\Psi_s$ maps the $i$-th vector in $\mathcal B_s$ to the $i$-th vector in $\mathcal B'_s$. Its matrix in the basis $\mathcal A$ is 
$$ P_s:=M_{\mathcal A\mathcal B'_s} M_{\mathcal B_s\mathcal A} =
\(|\ln s|^\dm M_{\mathcal A \mathcal B}\)^{-1}\left [|\ln s|^\dm M_{ \mathcal A\mathcal B}\right ]
 $$
and therefore $P_s$ tends to the identity matrix $I_3$ as $s\to 0$. 
Following the computation in the previous step, we define the fundamental period cell associated to the perturbed basis,
$$   C_s = \Phi_s^{-1}\left( \left[-\frac12,\frac12\right]^3\right)
    =\Psi_s^{-1}(\widehat C_s),   $$
and conclude that:
\be\label{usvs}  {G_s(u_s, A_s; \Omega)\over |\Omega|}=
   {G_s(u_s, A_s; C_s)\over |C_s|}\sim
   {G_s(v_s, B_s; \widehat C_s)\over |\widehat C_s|} \sim
    \frac{\lep^2}{2}\left( \alpha\|h\|_g + |h-H_{ex}|^2\right)
      .\ee
We  note for future reference that  
\be\label{psi} \Phi_s = |\ln s|^{1/2} \Phi_0\circ \Psi_s, \ee
where $\Phi_0$ is fixed and satisfies \eqref{vertical}, and $\Psi_s\to\text{Id}$ as $s\to 0$. 
\end{proof}

\begin{proof}[Step~3, Discretization] Next, we define a configuration $(\tilde u_{s,n}, A_s)$ for Lawrence--Doniach by discretization. We claim that there exists $t\in [0,s)$ such that, letting $\tilde u_{s,n}(x,y)=u_s(x,y,z_n+t)$, $\tau_t A_s(x,y,z)=A_s(x,y,z+t)$,  we have
\be\label{lg}   \mathcal{L}_{s,s}^\lambda (\tilde u_{n,s},\tau_t A_s)
    \le G_s(u_s, A_s). \ee
Note that this is {\em not} the desired final estimate, as we can only bound Lawrence--Doniach with $\eps=s$ in this way.  The full estimate will follow from adjusting $\tilde u_{n,s}$ in region $|\tilde u_{n,s}|\neq 1$.

To prove \eqref{lg} we use the following elementary observation:  for any continuous, $\mathcal{A}$-periodic real-valued function $f$, there exists $t\in [0,s)$ such that
\begin{equation}\label{trick}
  \int_\Omega f dx\, dy\, dz 
= \int_0^s \left[
    \sum_{n=1}^N \int_{P_n} f(x,y,z_n+\tau)\, dx\, dy\right] d\tau
= s \sum_{n=1}^N \int_{P_n} f(x,y,z_n+t)\, dx\, dy.
\end{equation}
The first equality comes from Fubini's theorem, and the second from the mean value theorem for integrals.  We apply the above with
$$  f= \frac12 |(d'-i A'_s)u_s|^2 + {1\over 4s^2}(|u_s|^2-1)^2,  $$
obtaining $t\in [0,s)$ for which \eqref{trick} holds.  We then define
$\tilde u_{s,n}(x,y)=u_s(x,y,z_n+t)$,  and conclude that 
\begin{multline}\label{trick1}
\sum_{n=1}^N \int_{P_n} \frac12 |(d'-i A'_s(x,y,z_n+t))\tilde u_{s,n}|^2 + {1\over 4s^2}(|\tilde u_{s,n}|^2-1)^2 \, dx\, dy \\
= \int_\Omega \frac12 |(d'-i \tau_tA'_s)u_s|^2 + {1\over 4s^2}(|u_s|^2-1)^2 dx\, dy\, dz.
\end{multline}
For the coupling term, we first calculate, 
\begin{align*}
{\tilde u_{s,n} -\tilde u_{s,n-1}e^{i\int_{z_{n-1}}^{z_{n}} \tau_t A_{s,z}(x,y,z')dz'}\over
   s} &=
{1\over s}\int_{z_{n-1}+t}^{z_n+t} \partial_z \left(
     u_s(x,y,z)e^{i\int_z^{z_n+t} A_{s,z}(x,y,z')dz'}\right) dz \\
     &= 
     {1\over s}\int_{z_{n-1}+t}^{z_n+t}
        (\partial_z u_s -i A_{s,z} v_s) e^{i\int_z^{z_n+t}A_{s,z}(x,y,z')dz'} dz.
\end{align*}
Then, applying Cauchy-Schwartz, summing over $n$, and using periodicity we find,
\be\label{trick2}
s \sum_{n=1}^N \int_{P_n} {1\over s^2}
   \left| \tilde u_{s,n} -\tilde u_{s,n-1}e^{i\int_{z_{n-1}}^{z_{n}}\tau_t A_{s,z}(x,y,z')dz'}
   \right|^2 dx\, dy \le 
   \int_\Omega |\partial_z u_s -i A_{s,z} u_s|^2.
\ee

Finally, since $A_s$ is $\mathcal{A}$-periodic, 
$$  \int_\Omega |d (\tau_tA_s) - h_{ex}^\eps|^2 =
\int_\Omega |d B_s - h_{ex}^\eps|^2 .
$$
Putting this together with \eqref{trick1} and \eqref{trick2} we obtain \eqref{lg}.
\end{proof}

\begin{proof}[Step~4, Conclusion] We may now prove Proposition~\ref{ubprop}.  For simplicity we denote from now on $\tau_t A_s$ by $A_s$. We modifiy $\tilde u_{s,n}$ in the set where $|\tilde u_{s,n}|\neq 1$. We have 
$$\{|\tilde u_{s,n}|\neq 1\} = P_n\cap {\Phi_s}^{-1}\(\cup_{p\in\ZZ^2} B(p,2s)\times\RR\).$$
From \eqref{vertical} and \eqref{psi}, the angle between $P_n$ and the recipocal image by $\Phi_s$ of $\vec e_3$ is bounded away from $0$, and therefore the intersection $E_{p,s,n}$ of each elliptical cylinder ${\Phi_s}^{-1}\(B(p,2s)\times\RR\)$ with $P_n$ is an ellipse whose aspect ratio depends only on $s$ and is bounded uniformly with respect to $s$. Still from \eqref{psi}, the diameter of the ellipse is bounded above by   $Cs|\ln s|^{-1/2}$. For small enough $s$, we may therefore include each ellipse in  the disk $D_{s,p}^n$ in the plane $P_n$ which has radius $s$ and  same center.  We let $T_{s,n} = \{D_{2s,p}^n\}_{p\in\ZZ^2}$ ($n=1\dots,N$) be the collection of disks with the same centers but twice the radius.

Outside these disks, $|\tilde u_{s,n}|=1$, and hence \eqref{lg}, \eqref{usvs}, \eqref{vs} allow to bound the energy in this region:
\begin{align}  \nnn
& s\sum_{n=1}^N 
 \int_{P_n\setminus T_{s,n}}\left\{\frac12 |(d' - iA'_s(x,y,z_n))\tilde u_n|^2
     + {1\over4\eps^2}(|\tilde u_{s,n}|^2-1)^2\right\} dx \, dy \\
     \nnn
     &\qquad + s\sum_{n=1}^N  {1\over 2\lambda s^2}
     \int_{P_n\setminus [T_{s,n}\cup T_{s,n-1}]}
      |\tilde u_{s,n} 
        -\tilde u_{s,n-1}
           e^{i\int_{z_{n-1}}^{z_n} A_{s,z}(x,y,z')dz'}|^2 dx\, dy
      \\
      \nnn
   &\qquad + \frac12\int_\Omega |dA_s - h_{ex}^\eps|^2\\
   \nnn
   &\qquad\qquad\le \mathcal{L}_{s,s}^\lambda(\tilde u_{s,n}, A_s)\\
   \label{outside}
     &\qquad\qquad\le G_s(u_s,A_s)\le 
\frac12\left( \alpha\| h\|_g + | h - H_{ex}|^2\right)|\ln\eps|^2 + o(|\ln\eps|^2)
\end{align}

To complete the upper bound we modify the functions $\tilde u_{s,n}$ inside each disk $D_{2s}$ to incorporate a two-dimensional vortex with core size $\eps$, and we call the result of this modification $u_{s,n}$. 

First, we define $u_{s,n}$ in $\omega_s:=D_{2s}\sm D_s$ to be of modulus $1$, and such that it  interpolates between $\tilde u_{s,n}$  outside $D_{2s}$ and $x\to (x-c)/|x-c|$ on $\partial D_s$, where  $c$ denotes the center of the $D_s$. From \eqref{defphi} and the discussion following it, near $c = \Phi_s^{-1}(p,z)$ we have, denoting by $\pi$ the projection on the plane $xy$ and $X = \pi\Phi_s(x)$, 
\be\label{forme}\frac{\tilde u_{s,n}(x)}{|\tilde u_{s,n}(x)|} = \frac{X-p}{|X-p|} e^{i\xi_p(X)}\ee
hence the gradient of $\tilde u_{s,n}$ is bounded by $C/s$ in $\omega_s$, and so is the gradient of $x\to (x-c)/|x-c|$. Since $A_s$ is bounded by $C|\ln s|$ we find that  
$$\int_{\omega_s}  |(d-iA'_s) u_{s,n}|^2 \le \frac C{s^2} |\omega_s|\le C.$$

Second, inside $D_s$, we let 
$$u_{s,n}(x) = \rho_\eps(|x-c|)\frac{x-c}{|x-c|},$$
where  $\rho_\eps=0$ on $[0,\eps]$, where $\rho_\eps(t) = (t-\eps)/\eps$ on $[\eps,2\eps]$ and  $\rho_\eps =1$ on $[2\eps,+\infty]$. We have
\begin{align}\label{machin}
&\int_{D_s} \frac12|(d' - iA'_s)u_{s,n}|^2 +{1\over 4\eps^2} (|u_{s,n}|^2-1)^2 \\ \nnn
&\qquad\le  (1+ |\ln s|^{-1}) \int_{D_s}
   \frac12 \left[|d (\rho_\eps e^{i\theta})|^2 + {1\over 4\eps^2} (\rho_\eps^2-1)^2\right]
 + (1+|\ln s|) \int_{D_s} \frac12 |A_s|^2 \\  \nnn
&\qquad
   \le \pi\ln \left({ s \over \eps}\right)[1+ o(1)],
\end{align}
where we have used $|x+y|^2\le (1+a)|x|^2 + (1+a^{-1})|y|^2$ (with $a=|\ln s|^{-1}$) in the second line.

It remains to count the number of disks in  $T_{s,n}$. From \eqref{forme} each disk in $T_{s,n}$ contributes a degree $1$ to the map $u_{s,n}$. On the other hand, since the configuration $(u_{s,n},A_s)$ is Floquet-periodic with respect to $(\vec v_1,\vec v_2)$, the total degree in each $P_n$ is calculated in terms of the magnetic flux,
\begin{align*}
  |\deg(u_{s,n},\partial P_n)| &= \left|{1\over 2\pi} \int_{P_n} dA_s \right|\\
    &= {|\ln s|\over 2\pi} \frac h\alpha\(\vec v_1,\vec v_2\) =  |\ln \eps|{|P_n|\over 2\pi} h\(\vec e_1, \vec e_2\) = |\ln \eps|{|P_n|\over 2\pi}h_3.
\end{align*}
Summing over the planes we obtain, since $s\sum_n|P_n| = |\Omega|$ and in view of \eqref{machin}
\begin{multline} \label{ellipses}
s\sum_{n=1}^N 
 \int_{T_{s,n}\subset P_n}\left\{ |(d' - iA'_s(x,y,z_n))u_{s,n}|^2
     + {1\over\eps^2}(|u_{s,n}|^2-1)^2\right\} dx \, dy 
   \\  \le {1-\alpha\over 2} |h_3|\, |\Omega|\, \lep^2 (1+o(1)).
\end{multline}
It remains to estimate the contribution of the excluded ellipses ${T}_{s,n}$ to the coupling term in the energy.  Since $|u_{s,n}|\le 1$, and from the above the measure of ${T}_{s,n}$ is of order $s^2|\ln s|$ (the area of an ellipse times their number) and  we conclude that
$$  s\sum_{n=1}^N  {1\over 2\lambda s^2}
     \int_{[ T_{s,n}\cup T_{s,n-1}]}
      |u_{s,n} 
        -u_{s,n-1}
           e^{i\int_{z_{n-1}}^{z_n} A_z(x,y,z')dz'}|^2 dx\, dy
            \le C \lep.
$$
Combining with \eqref{outside} and \eqref{ellipses}, we obtain the desired conclusion and Proposition~\ref{ubprop} is proven.
\end{proof}

\section{Lower bound}\label{LBsection}
\setcounter{thm}{0}

Let $(u_n,A)$ denote minimizers of $\ELDL$ for $s=1/N\to 0$.
From Proposition~\ref{ubprop}, $\ELDL(u_n,A)\le C |\ln\eps|^2$, so we immediately deduce an estimate for the induced field $h=h_s=dA$
 in particular
$$   \| h\|_{L^2}\le C \lep.  $$
In particular, passing to a subsequence we may conclude
$$   {h\over \lep}\weak H \qquad\text{in $L^2(\Omega;\Lambda^2(\RR^3))$}.
$$

We prove the following lower bound on minimizers:
\begin{thm}\label{LBthm}   Let $(u_n, A)$ be minimizers of $\ELDL$, with external field $h^\eps_{ex}$ satisfying \eqref{truc}.  
Then, there exists a sequence $N\to\infty$, $s=L/N\to 0$, $\eps=s^{1/\alpha}\to 0$ and periodic $H\in L^2_{loc}$ such that
$$  {h\over\lep}\weak H \quad\text{in $L^2(\Omega;\RR^3)$},  $$
and 
$$  \liminf_{s\to 0 \atop \eps\to 0} 
      {\ELDL(u_n,A;\Omega)\over |\Omega| \lep^2}  
      \ge \frac12\left[ (1-\alpha) |\overline{H}_3| 
           + \alpha\| \overline{H}\|_g + | \overline{H} - H_{ex}|^2\right], $$
where $\overline{H}$ is the componentwise average of $H$ over $\Omega$, with components
$$  \overline{H_i}:=\dashint_\Omega H_i 
$$
\end{thm}

\subsection{The planes}

For $\omega\subset P_n$, define
\begin{gather*}
E_n(u_n,A; \omega) = \int_{\omega} \frac12 |(d'-iA'_n)u_n|^2
     + {1\over 4\eps^2} (|u_n|^2-1)^2, \\
  J_n(u_n,A;\omega) = {1\over s^2}\int_\omega
     |u_n -u_{n-1}e^{i\int_{z_{n-1}}^{z_n} A_z(x,y,z')dz'}|^2 dx\, dy, \\
 M_n(A;\omega) = {1\over 2s}\int_{z_{n-1}}^{z_n} \int_\omega |dA-\hev|^2\, dx\, dy\, dz,
\end{gather*}
so that
$$  \ELDL(u_n,A)= s\sum_{n=1}^N [E_n(u_n,A;P_n)+J_n(u_n,A;P_n)+ M_n (A;P_n)].  $$
We distinguish those planes for which the contribution to the energy is reasonably small:  define
\be\label{Ieps}
I_\eps:= \left\{
  n\in [1, N]: \  E_n(u_n,A;P_n)+ M_n(A;P_n) \le \lep^3\right\}.
\ee
Since
$$  C\lep^2 \ge s\sum_{n\notin I_\eps}  E_n(u_n,A;P_n)+ M_n(A;P_n)
\ge s \,\text{card}\,(I_\eps^c) \lep^3,  $$
we conclude that 
$$\text{card}\,(I_\eps^c)\le C/(s\lep).  $$

  We obtain the following 
version of Theorem~4.1 of \cite{SSbouquin}:

\begin{thm}\label{4.1}
For all sufficiently small $\eps>0$ and for every $k\in I_\eps $, there exists a finite collection 
$\mathcal B_k=\mathcal B(\eps)= \{B_{k,i}\}_{i=1,\dots, m_k}$ of
closed (three-dimensional) balls such that
\begin{enumerate}
\item Each $B_{k,i}$ has center $a_{k,i}\in P_k\cap \Omega$ and radius
$r_{k,i}$ with $\sum_i r_{k,i}=s$. 
\item Let $V_k=\mathcal B_k\cap P_k$.  Then,
$$  \{x\in P_k: \ \left| |u_k|-1\right| \ge \eps^{\alpha\over 4}\}\subset V_k.  $$
\item Set $d_{k,i}= \deg(u_k, \partial B_{k,i}\cap P_k)$, and 
$D_k=\sum_{i=1,\dots, m_k} |d_{k,i}|$.  Then, there exists a constant $C$ such that
\begin{gather}\label{LB1}
D_k \le C \lep^2, \quad\text{and} \\
\label{LB2}
\frac12 \int_{\mathcal B_k\cap P_k}
\left[ |(d' -i A'_k)u_k|^2 + {1\over 2\eps^2}(1-|u_k|^2)^2\right]
   + \frac{s}{2}\int_{\mathcal B_k} |h|^2 
   \ge \pi D_k  \ln \left[{s\over \eps}\right]\,\left( 1 - o(1)\right).
\end{gather}
\end{enumerate}
\end{thm}

\begin{proof}
This result follows almost directly from the basic two-dimensional vortex ball construction in Chapter~4 of \cite{SSbouquin}.  The difference is that the Lawrence--Doniach energy does not exactly reduce to the two-dimensional Ginzburg--Landau functional, because the energy due to the magnetic field is integrated over the domain $\Omega\subset\RR^3$, and there is no control over its trace onto the planes $P_n$.  Thus, we must modify the construction of \cite{SSbouquin} so as to integrate $h$ over three-dimensional balls rather than planar disks.  

The only real change in the construction is in Lemma~4.4 of \cite{SSbouquin}, in which a lower bound is obtained for the energy of an $S^1$-valued map $v$.  Following the notation of \cite{SSbouquin}, we replace the definition of $X$ by
$$   X= \int_{S_r^+} dA,  $$
with $S_r^+$ the upper hemisphere of radius $r$, centered on $P_n$, with exterior unit normal vector $\nu$.  In the proof of Lemma~4.4, this quantity is integrated over the planar disk $\Delta_r$, and by Stokes' Theorem the two are in fact equal.  Our aim is an estimate of $|dA|^2$ on the surface, and here there is a minor difference in formula (4.18).  By applying Stokes and Cauchy-Schwartz, we obtain an estimate of the same  form as (4.16),
\be\label{4.16}
   \frac12 \int_{\partial S_r^+=\partial \Delta_r} |(d'-iA') v|^2 
      + {\lambda\over 2}\int_{S_r^+} |dA|^2 \ge
            \pi {|\deg(\partial \Delta_r)|^2\over r} \left( { 1\over 1+ {r\over \lambda}}\right).
\ee
In Proposition~4.3, this estimate is then integrated in $r$.  We define
$$   \mathcal F (x,r) = \frac12 \int_{\Delta_r} |(d'-iA') v|^2
              + {\lambda\over 2}\int_{B_r^+} |dA|^2,  $$
where we recall that $\Delta_r$ is a (planar) disk in $ P_n$, and
$B_r^+$ is the upper half ball (three-dimensional.)  The proof of Proposition~4.3 continues as before; the result is an estimate for $S^1$-valued maps $v$ on the family of balls $\{\mathcal B(t)\}_{t\in\RR_+}$,
which replaces (4.15) of \cite{SSbouquin},
$$  \frac12 \int_{\Delta_r\setminus\omega} |(d'-iA')  v|^2
              + {r_1-r_0\over 2}\int_{B_r^+} |d A|^2
              \ge \pi |\deg(v,\partial \Delta_r)| \ln {r_1\over r_2},  $$
and a similar change in estimate (4.14).  The remainder of the construction continues as before, except the magnetic energy (denoted by $H(A,\Omega)$ in \cite{SSbouquin}) will be integrated over solid balls, and the coefficient of this term in the final estimate will be linear (and not quadratic) in the ball radius $r$.
\end{proof}

We may now state the contribution to the lower bound coming from the vortex balls in the superconducting planes:

\begin{prop}\label{biplan}
With notation as in Theorem~\ref{4.1},
\be\label{biplaneq}
\liminf_{\eps\to 0} 
 {s\sum_{n=0}^{N-1} E_n (u_n, A; {\mathcal B}_n)\over \lep^2}
    \ge \frac12 (1-\alpha) \left|  \int_\Omega H_3\right|   .
\ee
\end{prop}

\begin{proof}
We first note that
$$   \sum_{k\in I_\eps} s^2\int_{\mathcal B_k} |h|^2
        \le  s^2\int_\Omega |h|^2 \le s^2\lep^2,  $$
and so this term in the lower bound \eqref{LB2} may be neglected.

Using the Floquet periodicity (\ref{tH*}), we note first that for every $k\in I_\eps$, we have flux quantization in each plane given by
\be \label{dk}  D_k\ge
\left|\sum_{i=1}^{ m_k} d_{k,i}\right|
=|\deg( u_k; \partial P_k)|=
\left|{1\over 2\pi}\int_{ P_k} h\right|.
\ee
We now claim that 
$$   \int_{\{z=t\}\cap\Omega} h = \int_{ P_k} h $$
is constant in $t$, $k=1,\dots,N$.
Indeed, by the regularity theory for minimizers of the Lawrence--Doniach functional (see \cite{BK}), $h_3=h(e_1,e_2)$ is continuous in $\Omega$, and the 2-form $h$ is smooth in the regions $\Omega_n:=\Omega\cap\{z\in (z_{n-1},z_n)\}$.  Let $t\in [z_{n-1},z_n)$, and define the horizontal slab,
$S_t:= \Omega\cap\{z\in (t,z_n)\}$.  We then apply the divergence theorem to $h$ in $S_t$.  By periodicity of $h$, the contribution from the lateral edges of $S_t$ cancel, and hence,
$$   \int_{\{z=t\}\cap\Omega} h = \int_{ P_n} h.  $$
Applying the above with $t=z_{n-1}$, we conclude that the flux through each plane is identical, and in fact is constant in $t$, and the claim is verified.
In particular, \eqref{dk} implies that
$$  2\pi D_k \ge \left|\int_{ P_k} h\right| 
= {1\over L}\left|\int_0^L \int_{\{z=t\}\cap\Omega} h_3\right|
={1\over L}\left|\int_\Omega h_3\right|.  $$

In consequence,
\begin{align*}
s\sum_{n=0}^{N-1}  E_n (u_n, A; {\mathcal B}_n) + \frac{s^2}{2}\int_{\mathcal B_k} |h|^2
&\ge s\sum_{n\in I_\eps}  E_n (u_n, A; {\mathcal B}_n) + \frac{s^2}{2}\int_{\mathcal B_k} |h|^2\\
&\ge s\sum_{n\in I_\eps}\pi D_k  \ln \left[{s\over \eps}\right]\,\left( 1 - o(1)\right)\\
&\ge \frac12 s\sum_{n\in I_\eps} \left|\int_{ P_k} h \right|
\ln \left[{s\over \eps}\right]\,\left( 1 - o(1)\right)\\
&= \frac12 \left|\int_\Omega h_3\right| \ {s\over L} 
   \left(\text{card}\,I_\eps\right) \ln \left[{s\over \eps}\right]\,\left( 1 - o(1)\right) \\
&= \frac12 \left|\int_\Omega h_3\right|
  \ln \left[{s\over \eps}\right]\,\left( 1 - o(1)\right),
\end{align*}
since $\text{card}\,I_\eps=N(1-o(1))$ and $Ns=L$.
We divide by $\lep^2$, and use $s=\eps^\alpha$ to obtain the result.
\end{proof}

\subsection{A comparison functional}
Take minimizers $(u_n,A)$ of $\ELDL$ and introduce the interpolation
\be\label{Psidef}  \Psi(x,y,z):= \left[
        {z-z_{n-1}\over s} \, u_n(x,y) + {z_n-z\over s}\,
           u_{n-1}(x,y) \eia\right]\, e^{-i \int_z^{z_n} A_z(x,y,z')\, dz'}.
\ee
We observe that if
$\omega_j$, $j=1,2,3$, are the functions associated to $(u_n,\vec 
A)\in \HH$ (see \eqref{tH})
then $\Psi, A$ satisfy the Floquet conditions
\be\label{psibc}  \Psi(\vec{x}+\vec{v}_j) = \Psi(\vec{x})e^{i\omega_j}, 
   \quad A(\vec{x}+\vec{v}_j)= A(\vec{x}) + d \omega_j,\quad j=1,2,3,  
\ee
with the same $\omega_j$. 
As in \cite{ABS1} we connect the terms of the Lawrence--Doniach energy to the covariant derivatives of $\Psi$.  For convenience we define:
\be\label{tdef}
t=t(z) = {z-z_{n-1}\over s} \quad\mbox{ for $z_{n-1}< z\le z_n$.}
\ee
Note that $t(z)$ is $s$-periodic in $z$.

The connection between the terms in $\ELDL$ and the Ginzburg--Landau energy density of $\Psi$ are given in the following lemma.  Since the estimates are done in local coordinates and there are no metric considerations, it will be convenient to express the quantities as vectors rather than forms.

\begin{lem}\label{compare}
Assume $(u_n,A)$ are minimizers of $\ELD$, and $\Psi$ is defined as in (\ref{Psidef}).
Then, for any point $(x,y,z)\in\RR^3$ and for any $\gamma>0$ we have:
\begin{align}
\label{zderiv}
&|(\partial_z - iA_z(x,y,z))\Psi(x,y,z)|^2 = {1\over s^2}\left|u_n-u_{n-1}\eia\right|^2 \\
\label{GL1}
&|(\grad' - i \vec A')\Psi|^2 \le
      (1+\gamma) \biggl\{
           t(z) |(\grad' - iA'(x,y,z_n))u_n|^2  \\
           \nnn
           &\qquad\qquad
            + (1-t(z))
                |(\grad' - iA'(x,y,z_{n-1}))u_{n-1}|^2 
                \biggr\} + (1+\gamma^{-1})|g_n|^2, \\
\label{GL2}
&\left( 1- |\Psi|^2\right)^2 \le
      2\left[ t(z) \left( 1- |u_n|^2\right)^2  
         + (1-t(z)) \left( 1- |u_{n-1}|^2\right)^2\right] \\
         \nnn &\qquad\qquad\qquad\qquad
             +\frac14 \left|u_n - u_{n-1}
               e^{i\int_{z_{n-1}}^{z_n} A_z(x,y,z')dz'}\right|^2,
\end{align}
where $g_n(x,y,z)$ are functions with the property
$$  \sum_{n=1}^N \int_\Omega |g_n|^2 \le C s^2 |\ln s|^2.  $$
\end{lem}

\begin{proof}  Let $(x,y,z)\in\Omega$.  The identity \eqref{zderiv} follows from an explicit calculation, which is done in Lemma~4.2 of \cite{ABS1}.  From the same proof we also recall that
\begin{align*}
\left( |\Psi|^2-1\right)^2  &=
 \left(   2t(1-t)\Re\{u_n^*u_{n-1}\eia\}
           + t^2 |u_n|^2 + (1-t)^2 |u_{n-1}|^2 -1 \right)^2 \\
&=  \left( -t(1-t)\left|  u_n-u_{n-1}\eia\right|^2 +
        t(|u_n|^2-1) + (1-t)(|u_{n-1}|^2-1)\right)^2  \\
&\le 2 t^2(1-t)^2 \left|  u_n-u_{n-1}\eia\right|^4 + 2\left[
    t(|u_n|^2-1) + (1-t)(|u_{n-1}|^2-1)\right]^2  \\
&\le \frac14 \left|  u_n-u_{n-1}\eia\right|^2 + 2\left[
    t(|u_n|^2-1)^2 + (1-t)(|u_{n-1}|^2-1)^2\right] ,
\end{align*}
by convexity in the last line, and using $|u_n|\le 1$, valid for minimizers.  Hence \eqref{GL2} holds.

Next, we calculate
\begin{align}\label{id1}
&(\partial_x-iA_x(x,y,z))\Psi = \\ 
\nnn
 & \qquad \bigg\{ t (\partial_x-iA_x(x,y,z_n))u_n    + (1-t)\eia(\partial_x-iA_x(x,y,z_{n-1}))u_{n-1}\bigg\}\times \\
 \nnn &\qquad\qquad
   \times\eiaz  + g_{n},
 \end{align}
where $g_n = g_{n,1}+g_{n,2}$, 
\begin{gather*}
g_{n,1}=  itu_n  \eiaz \int_z^{z_n} h_y\, dz'\\
g_{n,2}=   - i(1-t)u_{n-1}e^{i\int_{z_{n-1}}^z A_z(x,y,z')\, dz'}
\int_{z_{n-1}}^ z  h_y \, dz'  
\end{gather*}
We then estimate the error terms,
\begin{align*}
\sum_{n=1}^N \int_{\Om_n} |g_{n,1}|^2 \, dz\, dx\, dy 
&\le  \sum_{n=1}^N \int_{\Om_n} 
  \left| \int_z^{z_n} h_y\, dz'\right|^2 dz\, dx\, dy \\
  &\le s \sum_{n=1}^N \int_{\Om_n} \int_z^{z_n} h_y^2
   \, dz' \, dz\, dx\, dy \\
 &\le s^2 \sum_{n=1}^N \int_{ P_n}\int_{z_{n-1}}^{z_n} h_y^2
   \, dz' \, dx\, dy \\
 &= s^2 \|h_y\|_{L^2(\Omega)}^2 
 \le s^2\  \lep^2= o(1).
\end{align*}
A similar estimate holds for $g_{n,2}$.  Expanding the identity \eqref{id1},  we conclude that, for any
$(x,y,z)\in\Om$ and 
$\gamma>0$,
\begin{align*}
&|(\partial_x -i A_x(x,y,z))\Psi|^2 \\
 &\quad \le (1+\gamma) \bigg|
     t (\partial_x-iA_x(x,y,z_n))u_n  
      + (1-t)\eia(\partial_x-iA_x(x,y,z_{n-1}))u_{n-1}\bigg|^2 \\ &\qquad
      +
       (1+\gamma^{-1})
        |g_n|^2 \\
&\quad\le 
  (1+\gamma) \bigg[  t |(\partial_x-iA_x(x,y,z_n))u_n|^2
      + (1-t)|(\partial_x-iA_x(x,y,z_{n-1}))u_{n-1}|^2 \bigg] \\ &\qquad
      +
       (1+\gamma^{-1})
       |g_n|^2.
\end{align*}
The analogous estimate for $(\partial_y-iA_y(x,y,z_n))u_n$ is proved in exactly the same way, with $-h_x$ playing the role of $h_y$.  These two bounds yield \eqref{GL1}.
\end{proof}
\medskip

We now use the interpolated function $\Psi$ to bound from below the energy in $\Omega\setminus \cup_n \mathcal B_n$, the region exterior to the balls constructed in 
Theorem~\ref{4.1}.  Since the value of $\Psi(x,y,z)$ in each gap $z_{n-1}<z<z_n$ depends on the values $u_n(x,y)$ and $u_{n-1}(x,y)$ in the two adjacent planes, excising a disk in the plane $P_n$ precludes calculating the value of $\Psi$ in a cylindrical region extending through {\em both} adjacent gaps, $z_{n-1}\le z\le z_{n+1}$.
Thus, in our lower bound we must consider the smaller domain
\begin{equation}\label{Omtilde}
\tilde\Om:=\Om\setminus  \mathcal C, 
\end{equation}
where  
$$  \mathcal C:= \bigcup_{n} \{ (x,y,z): \ z_{n-1}\le z\le z_{n+1} \ \text{and} \ 
(x,y,z_n)\in\mathcal B_n, \ 
   \}.
$$

We recall that the anisotropy is represented by the metric tensor $g=\text{diag}\,(1,1,\lambda^2)$, with
$|\omega|^2_g=\sum_{j,k} g^{jk}(\omega_j,\omega_k)$ for the 1-form $\omega = \omega_j dx^j$.

As in \cite{ABS1} we define (for $\omega\subset\RR^3$) the functional
\begin{equation}\label{MM} \left.\begin{gathered}
  \MM(\Psi,A,\omega) = \int_\omega m_s(\Psi, A),  \\
 m_s(\Psi,A):=\frac12 \biggl[ |\Psi|^2 |d_A\Psi|_g^2 + 
       (1-|\Psi|^2) \left( \frac12|\nabla |\Psi| |_g^2 + {2\beta\over s^2} (1-|\Psi|^2)^2\right)  \biggr] 
     \end{gathered}\right\}
\end{equation}  
with $\beta=\min\{1,\lambda^{-2}\}$.
We note that $\MM$ is very similar to the classical Ginzburg--Landau energy, but the role of the Ginzburg--Landau parameter is played by $2s^{-1}$.
Following the same steps as Lemma~4.4 of \cite{ABS1} we conclude that
\begin{prop}\label{micmac}
Let $(u_n,A)$ be minimizers of $\ELDL$, and $\Psi$ as in \eqref{Psidef}.
Then,
\begin{equation}\label{LBMM1}
  \ELDL(u_n,A)\ge \MM(\Psi,A; \Omega) 
   +\frac12 \int_\Omega |h-\he|^2 + o(1),  
\end{equation}
and (with $\tilde\Om$ as in \eqref{Omtilde},)
\begin{align} \label{LBMM2}
\MM(\Psi, A; \tilde\Omega) &\le
s\sum_{n=1}^{N(s)}
     \int_{ P_n\setminus\mathcal B_n} \left[  \frac12
            \left|\left( d' -i 
A'_n\right)u_n  \right|^2
+ \frac{1}{4\eps^2} (|u_n|^2-1)^2  \right]\, dx\, dy
          \\   \nnn
  &\qquad    + s\sum_{n=1}^N
                    \int_{ P_n}  {1\over 2\lambda^2 s^2}\left| u_n -
                      u_{n-1}\eia
                           \right|^2
                   \, dx\, dy + o(1) .
   \end{align}
\end{prop}


\begin{proof}   Let $(x,y,z)$ be chosen with $z_{n-1}<z\le z_n$.
From \eqref{zderiv} and \eqref{GL1}, with $\gamma=|\ln\eps|^{-4}$, we have
\begin{align*}
\frac14 \bigl| \nabla |\Psi|\bigr|^2_g 
&\le    \frac14 \bigl| \nabla_A \Psi\bigr|^2_g \\
&\le 
\frac14(1+\gamma) \bigl\{
           t(z) |(\grad' - iA'(x,y,z_n))u_n|^2   + (1-t(z))
                |(\grad' - iA'(x,y,z_{n-1}))u_{n-1}|^2 
                \bigr\} \\
  &\qquad + {1\over 4\lambda^2 s^2}
   \left|u_n-u_{n-1}\eia\right|^2
  + \frac14(1+\gamma^{-1})|g_n|^2 \\
  &\le  \frac14\biggl\{
           t(z) |(\grad' - iA'(x,y,z_n))u_n|^2  + (1-t(z))
                |(\grad' - iA'(x,y,z_{n-1}))u_{n-1}|^2 
                \biggr\} \\
  &\qquad + {1\over 4\lambda^2 s^2}
   \left|u_n-u_{n-1}\eia\right|^2
  + |\tilde g_n|^2,
\end{align*}
with $\tilde g_n$ such that $\sum_{n=1}^N \int_\Om |\tilde g_n|^2 \to 0$.
Dividing \eqref{GL2} by $s^2\lambda^2$, and combining the resulting inequality with the above, we obtain:
\begin{align*}
&\frac14  \bigl| \nabla |\Psi|\bigr|^2_g 
      + {1\over s^2\lambda^2} (1-|\Psi|^2)^2   \\
      &\qquad\le
        \frac14\biggl\{
           t(z) |(\grad' - iA'(x,y,z_n))u_n|^2  + (1-t(z))
                |(\grad' - iA'(x,y,z_{n-1}))u_{n-1}|^2 
                \biggr\}   \\
      &\qquad + 
      {2\over s^2\lambda^2}\left[ t(z) \left( 1- |u_n|^2\right)^2  
         + (1-t(z)) \left( 1- |u_{n-1}|^2\right)^2\right] \\
      &\qquad  + {1\over 2\lambda^2 s^2}
   \left|u_n-u_{n-1}\eia\right|^2 + |\tilde g_n(x)|^2
\end{align*}
We may then conclude the pointwise bound,
\begin{align*}
m_s(\Psi,A)\le 
&\max\left\{ \frac12 |\nabla_A \Psi|_g, \  \frac14  \bigl| \nabla |\Psi|\bigr|^2_g 
      + {1\over s^2\lambda^2} (1-|\Psi|^2)^2\right\} \\
      &\le
       \frac12\biggl\{
           t(z) |(\grad' - iA'(x,y,z_n))u_n|^2  + (1-t(z))
                |(\grad' - iA'(x,y,z_{n-1}))u_{n-1}|^2 
                \biggr\}   \\
      &\qquad\qquad + 
      {1\over 4\eps^2}\left[ t(z) \left( 1- |u_n|^2\right)^2  
         + (1-t(z)) \left( 1- |u_{n-1}|^2\right)^2\right] \\
      &\qquad\qquad  + {1\over 2\lambda^2 s^2}
   \left|u_n-u_{n-1}\eia\right|^2 + 2|\tilde g_n(x)|^2,
\end{align*}
for all $(x,y,z)$ with $z_{n-1}<z\le z_n$, as long as $\eps^2<s^2\lambda^2/8$.

The lower bound \eqref{LBMM1} then follows by integrating over each gap and summing over $n$, using fact that $t(z)$ is $s$-periodic to combine adjacent terms in the sum and obtain $\ELDL(u_n,A)$ on the right-hand side.
For the bound \eqref{LBMM2} we must be more careful, since we cannot include contributions to the Lawrence--Doniach energy from the vortex balls in the planes.  Notice however that if we integrate only over $(x,y,z)\in \tilde\Omega=\Omega\setminus \mathcal{C}$, the vortex balls in the adjacent planes above and below are indeed excluded, and thus integrating and summing over the region $\tilde\Omega$ again yields the desired bound \eqref{LBMM2}.
\end{proof}

\subsection{Lower bound by slicing}

\begin{prop}\label{lbMM}  Let $\Psi$ be defined as in \eqref{Psidef}, and $\MM$ as in \eqref{MM}.  Then,
$$  \lim_{s\to 0} {\MM(\Psi,A;\tilde\Omega)\over |\Omega| |\ln \eps|^2}
  \ge \frac{\alpha}2 \|\bar H\|_g,
$$
where the 2-form $\bar H$ is the component-wise average of the limiting field $H$.
\end{prop}

We derive the lower bound on $\MM(\Psi,A;\tilde\Omega)$ using the slicing method of \cite{ABS2}.  First, denote by $\widehat{\mathcal C}$ the periodic
extension (by integer multiples of $\vec{v_1},\vec{v_2},\vec{v_3}$) of the exceptional set $\mathcal C$  to $\RR^3$.
Then, fix a basis $(W_1,W_2,W_3)$ of $\RR^3$ consisting of integer linear combinations of the basis $(\vec{v_1},\vec{v_2},\vec{v_3})$ which generates the period domain $\Omega$.   Define the domain spanned by this new basis,
$$  Q:= \{ \alpha_1 W_1 +\alpha_2W_2 +\alpha_3 W_3: \ 0\le \alpha_j\le 1, \ j=1,2,3\},  $$
and the perforated domain,
$$  \tilde Q:= Q\setminus \widehat{\mathcal C}.  $$
It is then immediately true that:

\begin{lem}\label{newper}
$\Psi$ is Floquet-periodic with respect to the frame $(W_1,W_2,W_3)$; that is,
there exist $\xi_j\in H^2_{loc}(\RR^3)$, $j=1,2,3$ so that
$$  \Psi(\cdot + W_j)=\Psi(\cdot)e^{i\xi_j(\cdot)}, \qquad
A(\cdot + W_j)= A(\cdot) + d \xi_j(\cdot), \quad j=1,2,3 $$
holds in $Q$.
Moreover,  
\begin{gather}\label{energyper} 
 {1\over |Q|} \MM(\Psi,A;Q) = 
     {1\over |\Omega|} \MM(\Psi,A; \Omega), \quad\text{and}
     \\
      {1\over |Q|} \MM(\Psi,A;\tilde Q) = 
     {1\over |\Omega|} \MM(\Psi,A;\tilde \Omega). \label{energyper2}  
\end{gather}
\end{lem}
Note that \eqref{energyper} follows from the fact that the energy density of $\MM$ is periodic with respect to both domains $\Omega$ and $Q$.  To conclude \eqref{energyper2}, we write $\MM(\Psi,A;\tilde \Omega)=\int_\Omega m_s (\Psi, A) \chi_{\tilde\Omega}$, and note that (by definition of  $\widehat{\mathcal C}$) the density $m_s (\Psi, A) \chi_{\tilde\Omega}$ is also periodic.

We now foliate the domain $\tilde Q$ by planes 
$\Pi_t$ on which $\alpha_3=t$, $0\le t\le 1$.  Denote by $d\sigma_g$ the surface measure on the plane $\Pi_t$ in the anisotropic metric $g$.  The (Euclidean) volume measure $dV=dx\, dy\, dz$ on $Q$ is proportional to 
$d\sigma_g\, d\alpha_3$, with constant of proportionality
$dV(W_1,W_2,W_3)/d\sigma_g(W_1,W_2)$, that is,
\begin{equation}\label{measure}
dV = {|Q|\over |W_1\wedge W_2|_g} d\sigma_g\, d\alpha_3.
\end{equation}

The next step is to estimate the energy in each plane $\Pi_t$ using the vortex-ball construction.  In order to do this, we need verify two conditions:  first,
the planes $\Pi_t$ are perforated by the exceptional sets $\widehat{\mathcal C}$, 
and these holes must be accounted for by vortex balls.  To do this, the total radius of the excluded holes must be sufficiently small.  Secondly, the vortex ball construction demands that the (two-dimensional) energy in the plane $\Pi_t$ be not too large. 

\begin{lem}\label{Tlem}
There exists a measurable subset $T\subseteq [0,1]$ so that:
\begin{enumerate}
\item  for every $t\in T$, we have both
\begin{equation}\label{cond1}    \frac12\int_{\Pi_t} \biggl\{ |\Psi|^2 |d_A\Psi|_g^2 + 
     \frac12  (1-|\Psi|^2) \left( |\nabla |\Psi| |^2 + {4\beta\over s^2} (1-|\Psi|^2)^2\right)  \biggr\} \le s^{-\frac12} 
     \end{equation}
and the total radius of the exceptional set,
\begin{equation}\label{cond2}
\rho(\Pi_t\cap \widehat{\mathcal C})\le s|\ln s|^4.
\end{equation}
\item  The measure of the complement, $|T^c|\le C|\ln s|^{-2}$.
\end{enumerate}
\end{lem}
\begin{proof}
First, note that by \eqref{energyper} and Proposition~\ref{micmac}, we have
$\MM(\Psi,A;Q)\le C |\ln s|^2$, so by Fubini's Theorem the measure of the set $\tau_1$ of $t\in [0,1]$ for which \eqref{cond1} fails is of order $s^{1/2}|\ln s|^2$.  To estimate the measure of the set 
$\tau_2$ of $t\in [0,1]$ for which 
$\rho(\Pi_t\cap \widehat{\mathcal C})> s|\ln s|^4$, we list the cylinders in the exceptional set $\widehat{\mathcal C}=\{ C_k\}_{k=1,\dots, K}$, and note that there are at most $K=N\lep^2$ cylinders (one for each vortex ball), and each has radius at most $2s$.  Then, taking $\rho(C_k\cap\Pi_t)=0$ when $C_k\cap\Pi_t=\emptyset$, we have
$$   \int_0^1 \rho(\Pi_t\cap \widehat{\mathcal C}) \, dt =
              \sum_{k=1}^K \int_{C_k\cap\Pi_t\neq\emptyset} 
              \rho(C_k\cap\Pi_t)\, dt 
              \le C (2s)^2 \, N\lep^2 \le C' \, s (\ln s)^2.
$$
The set $\tau_2$ of $t$ for which $\rho(\Pi_t\cap \widehat{\mathcal C})> s|\ln s|^4$ is thus at most of measure $|\ln s|^{-2}$.  
We conclude that $T^c=\tau_1\cup\tau_2$ has measure $|T^c|=O(|\ln s|^{-2})$.
\end{proof}

For $t\in T$ we apply the vortex ball construction as in \cite{SSbouquin}.  We have:

\begin{prop}\label{sliceballs}  For every $t\in T$, there exists a finite collection of disjoint closed disks $\mathcal D_t =\{\Delta_j^t\}_{j\in J_t}$ lying in the plane $\Pi_t$ such that:
\begin{enumerate}
\item the sum of the radii 
$\rho(\mathcal D_t)= \sum_{j\in J_t} \rho(\Delta_j^t) = |\ln s|^{-10}:= r$;
\item $\{x\in \Pi_t: \ \bigl| |\Psi|^2-1\bigr|^2 \ge s^{1/8}\} \subset \mathcal D_t$;
\item $\widehat{\mathcal C}\cap \Pi_t\subset \mathcal D_t$;
\item  Let $d_j^t = \deg (\Psi; \partial \Delta_j^t)$ and 
$D_t = \sum_{j\in J_t} |d_i^t|$.  Then,
\begin{equation}\label{slicelb}
   \int_{\mathcal D_t \setminus \widehat{\mathcal C}}
 \left[  m_s(\Psi,A) + \frac12 r^2 |dA|_g^2  \right]
  d\sigma_g \ge \pi D_t |\ln s| (1- o(1)).
\end{equation}
\end{enumerate}
\end{prop}

The proof of Proposition~\ref{sliceballs} follows from Theorem~4.1 of \cite{SSbouquin}.  The primary difference in the two proofs that the exceptional set
$\widehat{\mathcal C}\cap\Pi_t$ must be included in the collection the initial balls
($\mathcal B_0$ in the notation of \cite{SSbouquin}.)  For $t\in T$, the total radius of the balls satisfies \eqref{cond2}, and so including these sets does not affect the remainder of the proof.  The density $m_s(\Psi,A)$ may be substituted for the usual Ginzburg--Landau energy density since $|\Psi|\ge 1-s^{1/8}$ outside of the initial balls ($\mathcal B_0$.)

\medskip

We note that by conclusion 1 of the Proposition~\ref{sliceballs} and the periodicity of $|\Psi|$ we may (if necessary) translate the period domain $Q$ in order that no vortex ball in $\Pi_t$ intersects the boundary.  As a consequence, we may define the degree of $\Psi$ along this boundary, $\deg(\Psi;\partial\Pi_t)$.  By Floquet periodicity and Stokes Theorem, the degree quantizes the flux of the magnetic field $h=dA$ through the surface $\Pi_t$:
\begin{equation}\label{quantize}
\int_{\Pi_t} h = 2\pi \deg(\Psi; \partial\Pi_t) \le 2\pi D_t.
\end{equation}
Thus, the lower bound \eqref{slicelb} implies
\begin{equation}\label{slice2}
 \int_{\mathcal D_t \setminus \widehat{\mathcal C}}
 \left[  m_s(\Psi,A) + \frac12 r^2 |dA|_g^2  \right]
  d\sigma_g \ge \frac12 \left(\int_{\Pi_t} h \right)  |\ln s| (1- o(1)).
\end{equation}

We now integrate the lower bound \eqref{slice2} over $t\in T$, using the identity \eqref{measure}:
\begin{align*}
& \int_{\tilde Q} 
  \left[   m_s(\Psi,A) + \frac12 r^2 |dA|_g^2 \right] dV
  \\
  &\qquad \ge \frac12 \int_{t\in T} 
    \int_{\mathcal D_t \setminus \widehat{\mathcal C}}
 \left[  m_s(\Psi,A) + r^2 |dA|_g^2  \right]
 {|Q|\over |W_1\wedge W_2|} d\sigma_g\, dt \\
 &\qquad \ge \frac12 {|Q|\over |W_1\wedge W_2|}\left(\int_T \int_{\Pi_t} h \right)  |\ln s| (1- o(1)) \\
 &\qquad = \frac12 {|Q|\over |W_1\wedge W_2|}\left(\int_0^1 \int_{\Pi_t} h \right)  |\ln s| (1- o(1)),
\end{align*}
since by Lemma~\ref{Tlem} we have
$$  \left|\int_{T^c}\int_{\Pi_t} h\right| \le |T^c|^{1/2} \|h\|_{L^2} \le O(1).
$$
Now, since we have
$$   \int_0^1 \int_{\Pi_t} h  = \int_0^1 \int_0^1\int_0^1 h(W_1,W_2)\, d\alpha_1 d\alpha_2\, dt = {1\over |Q|} \int_Q h(W_1,W_2) dV,  $$
the right-hand side simplifies to
$$  \frac12 \int_Q  {h(W_1,W_2)\over |W_1\wedge W_2|_g} dV.  $$
Finally, we remark that the term $r^2 |dA|_g^2$ above may be neglected.  Indeed,
$$  \int_{\tilde Q} r^2 |dA|_g^2 \le C r^2 [\ELDL(u_n,A) + |Q||\he|^2] 
\le o(1).  $$
Therefore, we rewrite the lower bound as:
$$  \MM(\Psi,A; \tilde Q) \ge 
     \frac12 \int_Q  {h(W_1,W_2)\over |W_1\wedge W_2|_g} dV \, 
           |\ln s| \, (1-o(1)).
$$
Writing $\bar h$ as the constant 2-form representing component-wise average of the 2-form $h$, we have ${\bar h\over \lep}\to \bar H$, the component-wise average of the limiting field.  Note that by periodicity, the averages are the same for any period domain, and so (by Lemma~\ref{newper}) we conclude
$$  \liminf_{s\to 0} {\MM(\Psi,A;\tilde\Omega)\over |\Omega| |\ln s\eps|^2}
= \liminf_{s\to 0} {\MM(\Psi,A;\tilde Q)\over |Q| |\ln \eps|^2}
  \ge \frac{\alpha}2 {\bar H(W_1,W_2)\over |W_1\wedge W_2|_g} 
  \ge \frac{\alpha}2 {\bar H(W_1,W_2)\over |W_1|_g|W_2|_g}.
$$
Taking the supremum over all such bases $(W_1,W_2,W_3)$ which are integer multiples of the periods we obtain the desired lower bound, and Proposition~\ref{lbMM} is complete.

\subsection{Finishing the proof of Theorem~\ref{bigthm}}

Putting together the results above, and using $s=\eps^\alpha$, we have:
$$  \liminf_{s\to 0 \atop \eps\to 0} 
      {\ELDL(u_n,A;\Omega)\over |\Omega| \lep^2}  
      \ge F(\overline{H})=\frac12\left[ (1-\alpha) |\overline{H}\cdot e_3| 
           + \alpha\| \overline{H}\|_g + | \overline{H} - H_{ex}|^2\right]. $$
Applying the upper bound from Proposition~\ref{ubprop} with $h=H_*$ the minimizer of $F(H)$ over $H\in\RR^3$ produces the complementary inequality,
$$  \limsup_{s\to 0 \atop \eps\to 0} 
      {\ELDL(u_n,A;\Omega)\over |\Omega| \lep^2}  
      \le F(H_*)\le F(\overline{H}). $$
In particular, $H(x)=\overline{H}=H_*$ is constant, and the normalized energy converges, proving \eqref{conv}.  This completes the proof of Theorem~\ref{bigthm}.  

\section{Minimizers of the limit energy}

In this section we use convex duality to derive a more geometrical characterization of the limiting normalized induced field, $H$.  In this setting, it will be convenient to represent the 2-form as a vector field, with the appropriate norms as defined in the Introduction.

\begin{prop}\label{dual}
Let
$$  F(H)= \frac12 |H-H_{ex}|^2 + {1-\alpha\over 2}|H_3|
                + {\alpha\over 2} \|H\|_g  $$
Then the Fenchel dual is
$$  F^*(H) = \begin{cases}
       \frac12 |H|^2 - \frac12 |H_{ex}|^2, & \text{if $H-H_{ex}\in K$,}\\
+\infty, &\text{if $H-H_{ex}\not\in K$,}
\end{cases}
$$
where $K$ is defined by
\be\label{claim1}   U\in K \iff \begin{cases}
\| U-{1-\alpha\over 2} e_3\|_{g^{-1}}\le {\alpha\over 2}, 
  &\text{if $U_3\ge {1-\alpha\over 2},$} \\
\| U+{1-\alpha\over 2} e_3\|_{g^{-1}}\le {\alpha\over 2}, 
  &\text{if $U_3\le -\left({1-\alpha\over 2}\right),$} \\
\| U'\|_{g^{-1}}\le {\alpha\over 2} ,
   &\text{if $-\left({1-\alpha\over 2}\right) \le U_3\le {1-\alpha\over 2},$}
\end{cases}
\ee
where we denote by $U'=(U_1,U_2,0)$ and $U_3=U\cdot e_3$.
Thus, if $H_*\in \RR^3$ minimizes $F$, then it also minimizes $F^*$, and $F(H_*)=-F^*(H_*)$.
\end{prop}

\begin{proof}
Let $U=H-H_{ex}$, equivalently, we minimize
$$  G(U) = \frac12 |U|^2 + \Psi(U),
\qquad  \Psi(U)= {1-\alpha\over 2}|U_3 + H_{ex,3}| +
                                 {\alpha\over 2} \|U+H_{ex}\|_g.  $$
By convex duality,
$$  \min_{H\in\RR^3} F(H) = \min_{U\in\RR^3} G(U)
= -\min_{U\in\RR^3} \frac12 |U|^2 + \Psi^*(-U),  $$
where $\Psi^*$ is the Fenchel dual of $\Psi$,
\begin{align*}
\Psi^*(U) &= \sup_{V\in\RR^3} \left[ U\cdot V - \Psi(V)\right] \\
&= - \inf_{V\in\RR^3} \left[ -U\cdot V +\Psi(V)\right] \\
&= -U\cdot H_{ex} - \inf_{V\in\RR^3} \left[
      -U\cdot V + {1-\alpha\over 2} |V\cdot e_3| +
         {\alpha\over 2} \| V\|_g \right].
\end{align*}
If there exists $V\in\RR^3$ with 
$-U\cdot V + {1-\alpha\over 2} |V\cdot e_3| +
         {\alpha\over 2} \| V\|_g<0$,
by homogeneity it follows that the infimum is $-\infty$, and then
$\Psi^*(U)=+\infty$.  On the other hand, if $U\in K$,
\be\label{good}
K:=\left\{U\in\RR^3: \ - U\cdot V + {1-\alpha\over 2} |V\cdot e_3| +
         {\alpha\over 2} \| V\|_g\ge 0,\ \text{for every $V\in\RR^3$,}
         \right\}
\ee
 then the infimum is attained with $V=0$, and 
$\Psi^*(U)= -U\cdot H_{ex}$.  

We claim that $K$ defined above coincides with \eqref{claim1}.
We divide into the three cases as above.

\noindent
{\bf Case $U_3\ge {1-\alpha\over 2}$:} \ 
First assume $U\in K$.  Let $V_0\in \RR^3$ be the vector with $\|V_0\|_g=1$ which attains the supremum,
$$   \| U-{1-\alpha\over 2} e_3\|_{g^{-1}}=
\sup_{\|V\|_g=1} (U- {1-\alpha\over 2} e_3)\cdot V =
   (U- {1-\alpha\over 2} e_3)\cdot V_0.  $$
Since $U_3-{1-\alpha\over 2}\ge 0$, clearly the sup is attained with $V_{03}=V_0\cdot e_3\ge 0$.  Thus,
$$  \| U-{1-\alpha\over 2} e_3\|_{g^{-1}} = (U- {1-\alpha\over 2} e_3)\cdot V_0 = U\cdot V_0 - {1-\alpha\over 2} |V_{03}| \le {\alpha\over 2},  $$
as desired.

Conversely, assume $\| U-{1-\alpha\over 2} e_3\|_{g^{-1}}\le {\alpha\over 2}$.  Then, for any $V\in \RR^3$, 
$$  {\alpha\over 2}\|V\|_g \ge (U-{1-\alpha\over 2}e_3)\cdot V
          \ge U\cdot V - {1-\alpha\over 2}|V_3|,  $$
 that is, $U\in K$.
 
 \medskip
 
\noindent
{\bf Case $U_3\le -\left({1-\alpha\over 2}\right)$:} \  
As above, assume $U\in K$, but now let $V_0$, $\|V_0\|_g=1$, attain
the supremum in
$$   \| U+{1-\alpha\over 2} e_3\|_{g^{-1}}=
\sup_{\|V\|_g=1} (U+ {1-\alpha\over 2} e_3)\cdot V =
   (U+ {1-\alpha\over 2} e_3)\cdot V_0.  $$
Now we must have $V_{03}=V_0\cdot e_3\le 0$, and hence
$$  \| U+{1-\alpha\over 2} e_3\|_{g^{-1}} = (U+ {1-\alpha\over 2} e_3)\cdot V_0 = U\cdot V_0 - {1-\alpha\over 2} |V_{03}| \le {\alpha\over 2},  $$
as desired.  The converse is proven exactly as in the previous case.

\medskip

\noindent

{\bf Case $-\left({1-\alpha\over 2}\right)\le U_3\le \left({1-\alpha\over 2}\right)$:} \
First assume $U\in K$, and let $V_0=V'_0=(V_{01},V_{02},0)$ with $\|V_0\|_g=1$ attain the supremum in
$$  \|U'\|_{g^{-1}} = \sup_{V'\in \RR^2\times\{0\}\atop \|V'\|_g=1}
       U'\cdot V' = U'\cdot V'_0.  $$
Then,
$$  \|U'\|_{g^{-1}} =  U'\cdot V'_0  = U\cdot V_0 - {1-\alpha\over 2}|V_{03}| \le {\alpha\over 2},  $$
as desired.  

Conversely, if $\| U'\|_{g^{-1}}\le {\alpha\over 2}$, for any $V\in\RR^3$
we have
\begin{align*}
{\alpha\over 2} \|V\|_g \ge U'\cdot V &= U\cdot V - U_3 V_3 \\
&\ge U\cdot V - |U_3|\, |V_3| \\
&\ge U\cdot V - {1-\alpha\over 2}|V_3|.
\end{align*}
This concludes the proof of the claim \eqref{claim1}.

We conclude that 
$$  \Psi^*(U) = \begin{cases}
-U\cdot H_{ex}, & \text{if $U\in K$,}\\
+\infty, &\text{if $U\not\in K$,}
\end{cases}
$$
and reverting to the original variable $H$ we conclude that
$$  F^*(H) = \begin{cases}
       \frac12 |H|^2 - \frac12 |H_{ex}|^2, & \text{if $H-H_{ex}\in K$,}\\
+\infty, &\text{if $H-H_{ex}\not\in K$.}
\end{cases}
$$
\end{proof}

\end{document}